\documentclass[a4paper,12pt,double-sided,single-lined]{article}
\usepackage[utf8]{inputenc}
\usepackage[backend=biber,style=alphabetic,sorting=ynt]{biblatex}
\usepackage[toc,page]{appendix}
\usepackage{amsmath}
\usepackage[english]{babel}
\usepackage{amssymb}
\usepackage{amsthm}
\usepackage{tikz-cd}
\usepackage{IEEEtrantools}
\usepackage{graphicx}
\usepackage{subfig}
\theoremstyle{definition} \newtheorem{definition}{Definition}
\theoremstyle{plain} \newtheorem{lemma}[definition]{Lemma}
\theoremstyle{plain} \newtheorem{corollary}[definition]{Corollary}
\theoremstyle{plain}
\theoremstyle{plain} \newtheorem{theorem}[definition]{Theorem}
\theoremstyle{remark} \newtheorem*{remark}{Remark}
    \theoremstyle{remark}

\title{Deformation Quantization on $\mathbb{R}^d$}
\author{Haiqi Wu\\
        Master of Mathematics in Mathematics\\
        Trinity Term 2021}\date{}
\usepackage[parfill]{parskip}
\begin{document}
\maketitle
\begin{abstract}
    \noindent The primary aim of this essay, drawn from the author’s MMath dissertation at Oxford, is to present and explain Kontsevich’s formality theorem.\\
    The first two sections introduce the main topic. Sections 3 and 4 discuss Hochschild (co)homology, DGLA, and Maurer-Cartan elements. HKR theorem is stated in section 4, while a proof is deferred to appendix.
    In Section 5 we talk about $L_{\infty}$ morphisms. The central result of section 5 gives a penultimate step for solving the problem. We prove the main theorem in section 6.
\end{abstract}
\newpage
\tableofcontents
\newpage
\section{Introduction}
This paper aims at understanding Kontsevich's deformation quantization. In its origin, quantization provides a link between classical mechanics and quantum mechanics. On the other hand, a deformation of an object $S$ is a family of objects $\{S_t\}$ depending on the parameter $t$ with $S_{t_0}=S$. In our situation, we are given a manifold $M$ with a Poisson bracket $\{-,-\}$ on $C^{\infty}(M)$, the goal is to find a deformation $*_{\hbar}$ of the standard product $\cdot$, depending smoothly on a parameter $\hbar$, such that $*_{\hbar}$ coincide with $\cdot$ when $\hbar=0$, and
\[\lim_{\hbar\rightarrow 0}\frac{f*_{\hbar}g-g*_{\hbar}f}{\hbar}=\{f,g\}\]
for any $f,g\in C^{\infty}(M)$. \\
In section 2, we will discuss the Moyal product as a first case of deformation quantization. The main theorem is stated, on which all the following discussion will center. \\
In Section 3 we develop the theory of Hochschild (co)homology. This (co)homology theory of algebra encodes huge amount of information about the algebra. It turns out to be a powerful tool along our way of investigation. We will define a graded Lie algebra structure on Hochschild cochains, and show the compatibility between such Lie structure and the differential. By the end of the section we will give an equation (called Maurer-Cartan) which describe all star products.\\
Section 4 is about differential graded Lie algebras (DGLAs). We have Hochschild cochain as a first example of DGLA. We will see another important example of DGLA, namely the polyvector fields ($T_{poly}$). Maurer-Cartan elements will be introduced. Maurer-Cartan elements in $T_{poly}$ and $D_{poly}$ correspond to Poisson bivector fields and star products respectively. We will discuss the HKR theorem, which gives a quasi-isomorphism from $T_{poly}$ to $D_{poly}$. However, it fails to preserve the Lie bracket.\\
Section 5 collects some facts about $L_{\infty}$ algebras and their morphisms. It follows that $L_{\infty}$ morphisms are those morphisms between DGLAs we want to work with. The key result is that, if an $L_{\infty}$ morphism is a quasi-isomorphism as a chain map, then it preserves Maurer-Cartan elements modulo gauge group action. \\
Section 6 consists of proof for the main theorem on $\mathbb{R}^n$. We will construct a morphism $U$ between $T_{poly}$ and $D_{poly}$. After verifying all conditions for $U$ to be an extension of HKR map and to be $L_{\infty}$, we finish the proof for the main theorem on $\mathbb{R}^d$.
\subsubsection*{Conventions.} Throughout the thesis we assume manifolds are over $\mathbb{R}$ unless explicitly stated. All algebras $A$  are over a field $k$ of characteristic $0$ (usually take $k=\mathbb{R}$).  Homomorphisms are denoted as $Hom_{\square}(-,-)$, while occasionally $Lin(-,-)$ are used for emphasizing linearity over a field $k$.
\newpage
\section{Deformation Quantization}
In this section we define Poisson structures, star products, and state the main theorem which is the theme of this paper.
\subsection{Poisson structures}
We first give some definitions about Poisson structures.
\begin{definition}
    A \textit{Poisson Algebra} is an associative algebra $A$ together with a bracket
    $\{-,-\}:A\times A\rightarrow A$ satisfying
    \begin{enumerate}
        \item \textit{skew-symmetry:} $\{f,g\}=-\{g,f\}$, $\forall f,g\in A$;
        \item \textit{Jacobi identity:} $\{f,\{g,h\}\}+\{g,\{h,f\}\}+\{h,\{f,g\}\}=0$, $\forall f,g,h\in A$;
        \item \textit{Leibniz rule:} $\{f,gh\}=g\{f,h\}+\{f,g\}h$, $\forall f,g,h\in A$.
    \end{enumerate}
\end{definition}
\begin{remark}
Points $1$ and $2$ in above definition says that $\{-,-\}$ turns $A$ into a \textit{Lie Algebra}. Point $3$ says that $\{f,-\}:A\rightarrow A$ is a derivation for $f\in A$.\footnote{Recall from differential geometry, a \textit{derivation} of $C^{\infty}(M)$ is an $\mathbb{R}-$linear map $ d:C^{\infty}(M)\rightarrow C^{\infty}(M)$ satisfying Leibniz rule. We will discuss derivations in more detail.}
\end{remark}
Let $M$ be a ($C^{\infty}$) manifold, and $C^{\infty}(M)$ the algebra of smooth functions on $M$ with the standard point-wise multiplication $\cdot$ . \begin{definition}
    A \textit{Poisson Manifold} is a pair $(M, \{-,-\})$ where $M$ is a manifold and $\{-,-\}$ is a Poisson bracket on $C^{\infty}(M)$.
\end{definition}
\begin{remark}
    By skew-symmetry and Leibniz rule, $\{-,-\}$ comes from a bivector field on $M$. \textit{i.e.} $\exists \Pi\in \Gamma(\bigwedge\nolimits^2TM)$ such that for any $f$, $g$ in $C^{\infty}(M)$,
    \begin{equation*}
     \{f,g\}=\Pi(df,dg).
    \end{equation*}
In local coordinates $\{\frac{\partial}{\partial x_i}\}$, 
we can express above as \[\{f,g\}=\sum_{i,j}\Pi^{ij}\frac{\partial f}{\partial x_i}\frac{\partial g}{\partial x_j}\] where $1\leq i,j\leq dim(M)$. Such bivector field $\Pi$ are said to be \textit{Poisson}. Notice that by skew-symmetry, $\Pi^{ij}=-\Pi^{ji}\in C^{\infty}(M)$.Hence we also write $(M, \Pi)$ for a Poisson manifold, or even $M$ by abuse of language.
\end{remark}
\begin{definition}
Suppose $(M,\{-,-\})$, $(N,\{-,-\}')$ are Poisson manifolds. A \textit{Poisson map} is a ($C^{\infty}$) morphism $\phi : M\rightarrow N$ so that for any $f,g\in C^{\infty}(N)$, $x\in M$,
\[\{f,g\}'(\phi(x))=\{f\circ\phi,g\circ\phi\}(x). \]
When $M=N$ and $\phi$ is diffeomorphic, we write $\{-,-\}\sim\{-,-\}'$.
\end{definition}
It follows that $\sim$ is an equivalent relation on the family of Poisson structures on $M$. \newline\newline
\textbf{Examples of Poisson structure:}
 \begin{enumerate}
 \item \textit{(Constant Poisson structure)} Consider $M=(\mathbb{R}^{n}, \Pi)$ where $\Pi^{ij}$'s are constants, $\Pi^{ij}=-\Pi^{ji}$ $(1\leq i,j\leq n)$, so \[\Pi=\sum_{1\leq i,j\leq n}\Pi^{ij}\frac{\partial}{\partial x_i}\otimes\frac{\partial}{\partial x_j}.\]
 The bracket $\{-,-\}$ so defined is indeed Poisson: \\
 it is skew-symmetry because $\Pi^{ij}=-\Pi^{ji}$; \\
 it satisfies Leibniz rule because
 \begin{align*}
     \{f,gh\}&=\sum_{i,j}\Pi^{ij}\frac{\partial f}{\partial x_i}\frac{\partial (gh)}{\partial x_j}\\
     &=\sum_{i,j}\Pi^{ij}\frac{\partial f}{\partial x_i}(h\frac{\partial g}{\partial x_j}+g\frac{\partial h}{\partial x_j})\\
     &= g\{f,h\}+\{f,g\}h;
 \end{align*}
 it satisfies Jacobi identity because
 \begin{align*}
     \{f,\{g,h\}\}& =\sum_{i,j}\Pi^{ij}\frac{\partial f}{\partial x_i}\frac{\partial}{\partial x_j}( \sum_{k,l}\Pi^{kl}\frac{\partial g}{\partial x_k}\frac{\partial h}{\partial x_l})\\
     & =\sum_{i,j}\Pi^{ij}\frac{\partial f}{\partial x_i}(\sum_{k,l} \Pi^{kl}\frac{\partial}{\partial x_j}(\frac{\partial g}{\partial x_k}\frac{\partial h}{\partial x_l}))\\
     & =\sum_{i,j,k,l}\Pi^{ij}\Pi^{kl}\frac{\partial f}{\partial x_i}(\frac{\partial^2 g}{\partial x_j\partial x_k}\frac{\partial h}{\partial x_l}+\frac{\partial g}{\partial x_k}\frac{\partial^2 h}{\partial x_j \partial x_l}),
 \end{align*}
 similarly
 \begin{align*}
 \{g,\{h,f\}\}= \sum_{i,j,k,l}\Pi^{ij}\Pi^{kl}\frac{\partial g}{\partial x_i}(\frac{\partial^2 h}{\partial x_j\partial x_k}\frac{\partial f}{\partial x_l}+\frac{\partial h}{\partial x_k}\frac{\partial^2 f}{\partial x_j \partial x_l}),
 \end{align*}
 \begin{align*}
 \{h,\{f,g\}\}=\sum_{i,j,k,l}\Pi^{ij}\Pi^{kl}\frac{\partial h}{\partial x_i}(\frac{\partial^2 f}{\partial x_j\partial x_k}\frac{\partial g}{\partial x_l}+\frac{\partial f}{\partial x_k}\frac{\partial^2 g}{\partial x_j \partial x_l}),
 \end{align*} and adding them together cancel out all terms.
 \item\textit{(Lie-Poisson structure)} Suppose $(L,[-,-])$ is a finite dimensional Lie algebra and let $L^*$ be its dual.  We can define a Poisson structure on $L^*$ as follow:\\
  $\bullet$ Assume first that $f,g\in Lin(L^*,\mathbb{R})\subseteq C^{\infty}(L^*)$, where $Lin(L^*, \mathbb{R})= L^{**}\cong L$ is the dual space of $L^*$, in which case $df=f$ and $dg=g$. \\
  Let $\{f,g\}:= [f,g]$ under identification $L^{**}\cong L$. So \[\{f,g\}:L^*\overset{C^{\infty}}{\longrightarrow} \mathbb{R}\in L^{**} \]
  is defined as
  \[\alpha\mapsto \alpha [f,g], \]
  for $\forall\alpha\in L^*$ ;\\
 $\bullet$ In general suppose $f,g\in C^{\infty}(L^*)$. For $\forall\alpha\in L^* $, $df|_{\alpha}$ can be identified with an element of $L$. We then define
 \[\{f,g\}(\alpha):=\alpha [df|_{\alpha},dg|_{\alpha}],\]
 for $\forall\alpha\in L^*$.\\
 To check that $\{-,-\}$ satisfies Leibniz rule, we verify:
 \begin{align*}
     \{f,gh\}(\alpha)&=\alpha[df|_{\alpha},d(gh)|_{\alpha}]=\alpha[df|_{\alpha},gdh|_{\alpha}+(dg)h|_{\alpha}]\\
     &=g(\alpha)\alpha[df,dg]+h(\alpha)\alpha[df,dh]=g(\alpha)\{f,h\}(\alpha)+ \{f,g\}(\alpha)h(\alpha)\\
     &=g\{f,h\}(\alpha)+\{f,g\}h(\alpha).
 \end{align*}
  To check that $(C^{\infty}(L^*),\{-,-\})$ is a Lie algebra, we first notice that this is indeed the case on $Lin(L^*,\mathbb{R})$ (since $Lin(L^*,\mathbb{R})\cong L$ as Lie algebras), then extend to the whole of $C^{\infty}(L)$ by the differential $d$.
 \end{enumerate}
\subsection{Star products}
Given a commutative algebra $(A, \cdot)$ over a field $k$ and a formal variable $\hbar$, the aim of deformation quantization is to describe a (formal) product on $A[[\hbar]]$ that '\textit{deforms}' the original product $\cdot$. Precisely,
\begin{definition}
A \textit{(formal) deformed product} (or a \textit{star product}) is an associative, $k[[\hbar]]$-linear and $\hbar$-adic continuous\footnote{Naively, this is to assign with $A[[\hbar]]$ a topology so that $*_{\hbar}$ varies continuously with $\hbar$.} product $*_{\hbar}$ on $A[[\hbar]]$ of formal power series, which, on $A$ can be written as
\begin{align}
    a*_h b=\sum_{i\geq 0}B_i(a,b)h^i,
\end{align}
with $B_0(a,b)=a\cdot b$ is the original product.
\end{definition}
\begin{remark}
If we write elements of $A[[\hbar]]$ as \[\sum_{i\geq 0}a_i\hbar^i,\] then the the product $*_{\hbar}:A[[\hbar]]\bigotimes_{k[[\hbar]]}A[[\hbar]]\rightarrow A[[\hbar]]$ maps
\[(\sum_{i\geq 0}a_i\hbar^i,\sum_{j\geq 0}b_j\hbar^j)\mapsto \sum_{l\geq 0}(\sum_{i+j+k=l}B_k(a_i,b_j))\hbar^l,\] by expanding linearly over $k[[\hbar]]$.
\end{remark}
Notice that if we take $\hbar:=0$, the star product coincide with the original commutative product $\cdot$.
The condition on associativity means that, for any $a,b,c\in A$
\begin{align*}
\sum_{i+j=k}B_i(B_j(a,b),c)=\sum_{i+j=k}B_i(a,B_j(b,c)),
\end{align*}
where $Bi:A\times A\rightarrow A$ for all $i$.\\
In particular, we have $B_1(ab,c)+B_1(a,b)c=B_1(a,bc)+aB_1(b,c)$. If we define $B_1^-(a,b):=(B_1(a,b)-B_1(b,a))/2$, then $B_1^-$ satisfies skew-symmetry. It can be checked that $B_1^-$ also satisfies Jacobi identity and Leibniz rule by expanding out all the terms. 
Hence $B_1^-$ is a Poisson bracket on $A$. In other words, if we have a deformed product $*_{\hbar}$ of $(A,\cdot)$, then
\[\{a,b\}:=\lim_{\hbar\rightarrow 0}\frac{a*_{\hbar}b-b*_{\hbar}a}{\hbar}\]
is a Poisson bracket on $A$.\\
Consider a Poisson manifold $M$. The space $A:=C^{\infty}(M)$ is equipped with the standard point-wise commutative product $\cdot$ and a Poisson bracket $\{-,-\}$. We are interested in a partial converse of above statement.
\begin{definition}
A \textit{deformation quantization} of Poisson manifold $M$ is a star product $*$ deforming $\cdot$ such that
\begin{enumerate}
    \item All $B_i$ in $(1)$ are bilinear bidifferential operators;
    \item $B_1^-$ coincides with $\{-,-\}$.
\end{enumerate}
\end{definition}
We can now state our main theorem
\begin{theorem}
\cite{Kontsevich_2003}
: Every Poisson manifold admits a deformation quantization.
\end{theorem}
Similar to Poisson brackets, there is an equivalent relation on the family of star products of $M$:
\begin{definition}
Two star products $*$ and $*'$ on $C^{\infty}(M)$ are said to be equivalent whenever there is a a $\mathbb{R}[[\hbar]]$-linear operator $D:C^{\infty}(M)[[\hbar]]\rightarrow C^{\infty}(M)[[\hbar]]$, which on $C^{\infty}(M)[[\hbar]]$ has the form
\[D(f)=f+\sum_{i\geq 1}D_if\hbar^i\] for $f\in C^{\infty}(M)$, such that
\[f*'g=D^{-1}(Df*Dg).\]
\end{definition}
Since $D=id+D_1\hbar+D_2\hbar^2+\cdots$, it starts with a unit thus is indeed invertible hence the equivalence above makes sense. In fact, they form a \textit{gauge group} acting on the family of star products.\\
\subsubsection{Example: Moyal product}
Consider the previous example of constant Poisson structure $(\mathbb{R}^n,\Pi)$. We define the star product by
\begin{align}
    f*g(x):=exp(\sum_{1\leq i,j\leq n}\hbar\Pi^{ij}\frac{\partial}{\partial x_i}\frac{\partial}{\partial y_j})(f(x)g(y))|_{y=x} & & f,g\in C^{\infty}(\mathbb{R}^n)
\end{align}
Then the star product so defined is associative:
\begin{align*}
(f*g)*h&=exp(\sum_{1\leq i,j\leq n}\hbar\Pi^{ij}\frac{\partial}{\partial x_i}\frac{\partial}{\partial z_j})((f*g)(x)h(z))|_{x=z}\\
       &=exp(\sum_{1\leq i,j\leq n}\hbar\Pi^{ij}(\frac{\partial}{\partial x_i}+\frac{\partial}{\partial y_i})\frac{\partial}{\partial z_j})exp(\sum_{1\leq k,l\leq n}\hbar \Pi^{kl}\frac{\partial}{\partial x_k}\frac{\partial}{\partial y_l})f(x)g(y)h(z)|_{x=y=z}\\
       &=exp(\sum_{1\leq i,j,k,l \leq n}\hbar(\Pi^{ij}\frac{\partial}{\partial x_i}\frac{\partial}{\partial z_j}+\Pi^{kl}\frac{\partial}{\partial x_k}\frac{\partial}{\partial y_l}+\Pi^{mn}\frac{\partial}{\partial y_m}\frac{\partial}{\partial z_n}))f(x)g(y)h(z)|_{x=y=z}\\
       &=f*(g*h)(x).
\end{align*}
The product $*$ is in fact a deformation quantization of $\{-,-\}$ defined by $\Pi$, since
\begin{align*}
(f*g-g*f)(x)&=\hbar\sum_{i,j}(\Pi^{ij}\frac{\partial f}{\partial x_i}\frac{\partial g}{\partial x_j}-\Pi^{ij}\frac{\partial g}{\partial x_i}\frac{\partial f}{\partial x_j})+O(\hbar^2) & & & & (\Pi^{ij}=-\Pi^{ji}) \\
            &=2\hbar\{f,g\}.
\end{align*}
 \subsection{Formal Poisson structures}
 Suppose $(A,\{-,-\})$ is a Poisson algebra and $\hbar$ is a formal variable. We can define a bracket $\{-,-\}_{\hbar}$ on $A[[\hbar]]$ to be a $k[[\hbar]]$-bilinear map
 \[\{-,-\}_{\hbar}:A[[\hbar]]\times A[[\hbar]]\rightarrow A[[\hbar]]\]
 so that $(A[[\hbar]],\{-,-\}_{\hbar})$ is a Poisson algebra.
\newline
For a Poisson manifold $(M,\Pi)$, we can define a bracket $\{-,-\}_{\hbar}$ on $C^{\infty}(M)[[\hbar]]$:
\[\{f,g\}_{\hbar}:=\sum_{l\geq 0}\hbar^l\sum_{i+j+k=l}\Pi_i(df_j,dg_k),\]
where $0\leq i,j,k\leq m$, $\Pi_i$ are bivector fields on $M$, $\Pi_0=\Pi$, and
\[f=\sum_{j\geq 0}f_j\hbar^j,g=\sum_{k\geq 0}g_k\hbar^k\in C^{\infty}(M)[[\hbar]].\]
\begin{definition}
Given $(M,\Pi)$ and $\Pi_{\hbar}$ as above, We call \[\Pi_{\hbar}:=\Pi_0+\Pi_1\hbar+\Pi_2\hbar^2+\cdots\] a \textit{formal Poisson structure} if $\{-,-\}_{\hbar}$ is a Poisson bracket on $C^{\infty}(M)[[\hbar]]$.\\
The gauge group in this case consists of formal power series of the form
\[\chi:=exp(\hbar X),\] called \textit{formal diffeomorphism},
where \[X=\sum_{i\geq 0}X_i\hbar^i\in\mathfrak{X}(M)[[\hbar]]\] is a formal power series over $\mathfrak{X}(M)$, called \textit{formal vector field}.\\
The group action is given by
\begin{align*}
\chi\cdot \Pi_{\hbar}&:=\chi\Pi_{\hbar}\chi^{-1}\\
&=exp(\hbar X)\circ\Pi_{\hbar}\circ(exp(-\hbar X)\otimes exp(-\hbar X)),
\end{align*}
where $\circ$ is the composition of functions.
\end{definition}
\begin{remark}
Let $A:=C^{\infty}(M)$, then a vector field $X_i$ in the definition of $X$ gives a linear function $A\rightarrow A$. Then $\chi=exp(\hbar X)$ givens a invertible $\mathbb{R}[[\hbar]]$-linear map $A[[\hbar]]\rightarrow A[[\hbar]]$.
We can represent the action of $\chi$ on $\Pi_{\hbar}$ via commutative diagram:
\[
\begin{tikzcd}[sep=1cm]
A[[\hbar]]\times A[[\hbar]] \arrow{r}{\chi\times \chi}\arrow[swap]{d}{\Pi_{\hbar}} &A[[\hbar]]\times A[[\hbar]]\arrow{d}{\chi\cdot\Pi_{\hbar}} \\
A[[\hbar]] \arrow{r}{\chi}  & A[[\hbar]] & .
\end{tikzcd}
\]
\end{remark}
Recall the \textit{Baker-Campbell-Hausdorff (BCH) formula} from a Lie Group course. For two vector fields $X_1,X_2$ on a manifold $M$, 
we have, by BCH,
\begin{align*}
  exp(X_1)exp(X_2)=exp(X_1+X_2+\frac{1}{2}[X_1,X_2]+\cdots).
\end{align*}
 This gives the group structure of $\{exp(\hbar X)|X\in\mathfrak{X}(M)[[\hbar]]\}$ by
 \begin{align*}
     exp(\hbar X)\cdot exp(\hbar Y):=exp(\hbar X+\hbar Y+\frac{1}{2}\hbar^2[X,Y]+\cdots).
 \end{align*}

We will follow Kontsevich's approach in his paper \cite{Kontsevich_2003}
to set up an identification between the set of formal Poisson structure modulo the gauge group action as above, and the set of star products modulo the gauge group action as mentioned in Definition $7$.
\newpage
\section{Hochschild (co)homology}
Hochschild (co)homology is a (co)homology theory on associative algebra, introduced by Gerhard Hochschild in 1945 in the classic paper \cite{10.2307/1969145}
. It turns out that Hochschild (co)homology is a powerful tool in analyzing certain differential graded algebraic structures. Moreover, the information of a \textit{deformed product} defined in previous section is encoded in Hochschild cochains, as we shall see. For more background information on the bar complex, please find appendix.\\
\subsection{Hochschild (co)chain complexes}
In the following we denote $A^e:= A\otimes_kA^{op}$ for an associative $k$-algebra $A$ as the \textit{enveloping algebra} of $A$. We regard, for any $n\in \mathbb{N}$, $A^{\otimes n}$ (tensoring over $k$) as an $A$-bimodule (or equivalently left $A^e$-module) by \[(a\otimes b)\cdot (a_1\otimes \cdots\otimes a_n):=aa_1\otimes\cdots\otimes a_nb.\]
\begin{definition}
Let $A$ be a $k$-algebra. The \textit{bar complex} $C_{\bullet}^{bar}(A)$ of $A$ is the chain complex
\[
\begin{tikzcd}
\cdots \overset{d_2}{\longrightarrow} A\otimes A \otimes A \overset{d_1}{\longrightarrow}  A \otimes A  \overset{d_0}{\longrightarrow} A \rightarrow 0,
\end{tikzcd}
\]
of $A$-bimodules (notice that $C_n^{bar}=A^{\otimes (n+2)}$), with the differentials $d_n:C^{bar}_n\rightarrow C^{bar}_{n-1}$ is defined by
\begin{equation}
    a_0\otimes \cdots \otimes a_{n+1}\mapsto \sum^n_{i=0}(-1)^ia_0\otimes \cdots \otimes a_{i-1}\otimes a_ia_{i+1}\otimes a_{i+2}\otimes \cdots \otimes a_{n+1}.
\end{equation}
In particular, $d_0$ maps $a_0\otimes a_1 $ to $a_0a_1\in A$.
\end{definition}
\begin{remark}
We can represent elements $ a_0\otimes \cdots \otimes a_{n+1}$ as $a_0[a_1|\cdots|a_n]a_{n+1}$, justifying the term \textit{bar complex}.
\end{remark}
Such $d_n$'s indeed makes $C^{bar}_{\bullet}$ into a chain complex, furthermore, $C_{\bullet}^{bar}$ gives a \textit{free resolution} of $A$ as an $A$-bimodule with \textit{augmentation map} $d_0$. We defer the detailed discussion into appendix.  \\
\newline
Suppose $M$ is a bimodule over a $k$-algebra $A$.
\begin{definition}
The \textit{Hochschild chain complex} $C_{\bullet}(A,M)$ is the complex of $k$-modules $C^{bar}_{\bullet}\otimes_{A^e}M$, with differentials $\delta_n=d_n\otimes id_M$.\\
The $n^{th}$-homology  group of $C_{\bullet}(A,M)$  
is denoted by $HH_n(A,M)$. When $M=A$ we simply write $HH_n(A)$.
\end{definition}
\begin{remark}
When $M=A^e$, $C_{\bullet}(A,M)=C_{\bullet}(A,A^e)=C^{bar}_{\bullet}(A)$.
\end{remark}
\begin{definition}
The \textit{Hochschild cochain complex} $C^{\bullet}(A,M)$ is the complex $Hom_{A^e}(C^{bar}_{\bullet},M)$, regarded as $k$-modules as well, with differential $\delta^n=Hom(d_n,id_M)$ (i.e. $f\overset{\delta^n}{\mapsto} f\circ d_{n+1}$).\\
The $n^{th}$-cohomology group of $C^n(A,M)$ is denoted by $HH^n(A,M)$. When $M=A$ we write $HH^n(A)$.
\end{definition}
Using (3) we can write out the differentials of a cochain explicitly. Suppose $f\in C^{n-1}(A,M)=Hom_{A^e}(A^{\otimes n+1},M)$, then $\delta^{n-1}(f)\in C^n(A,M)$ is unwrapped as:
\begin{align*}
    a_0\otimes\cdots a_{n+1}&\overset{d_n}{\longmapsto}\sum_{i=0}^n(-1)^ia_0\otimes\cdots a_{i-1}\otimes a_ia_{i+1}\otimes a_{i+2}\otimes\cdots\otimes a_{n+1}\\
    &=a_0a_1\otimes a_2\otimes\cdots\otimes a_n+(-1)^na_0\otimes\cdots\otimes a_{n-1}\otimes a_na_{n+1}\\
    &+\sum_{i=1}^{n-1}(-1)^ia_0\otimes\cdots a_{i-1}\otimes a_ia_{i+1}\otimes a_{i+2}\otimes\cdots\otimes a_{n+1}\\
    &=(a_0\otimes 1)(a_1\otimes\cdots\otimes a_{n+1})+(-1)^n(1\otimes a_{n+1})(a_0\otimes\cdots\otimes a_n)\\
    &+\sum_{i=1}^{n-1}(-1)^ia_0\otimes\cdots a_{i-1}\otimes a_ia_{i+1}\otimes a_{i+2}\otimes\cdots\otimes a_{n+1}\\
    &\overset{f}{\longmapsto}(a_0\otimes 1)f(a_1\otimes\cdots\otimes a_{n+1})+(-1)^n(1\otimes a_{n+1})f(a_0\otimes\cdots\otimes a_n)\\
    &+\sum_{i=1}^{n-1}(-1)^if(a_0\otimes\cdots a_{i-1}\otimes a_ia_{i+1}\otimes a_{i+2}\otimes\cdots\otimes a_{n+1})\\
    &=a_0f(a_1\otimes\cdots\otimes a_{n+1})+(-1)^nf(a_0\otimes\cdots\otimes a_n)a_{n+1}\\
    &+\sum_{i=1}^{n-1}(-1)^if(a_0\otimes\cdots a_{i-1}\otimes a_ia_{i+1}\otimes a_{i+2}\otimes\cdots\otimes a_{n+1}).\\
\end{align*}
The adjunction
\[-\otimes_k A^e: \text{ $k$-Mod}\rightleftharpoons \text{$A^e$-Mod }:res \]
gives an equivalence between the category of $k$-modules and the category of $A^e$-modules. Therefore, we have an isomorphism of $k$-modules
\[C^n(A,M)\overset{\cong}{\longrightarrow}Hom_k(A^{\otimes n},M)\]
given by
\[f\mapsto [a_1\otimes\cdots\otimes a_n\mapsto f(1\otimes a_1\otimes\cdots\otimes a_n\otimes 1)],\]
with inverse
\[g\mapsto [a_0\otimes\cdots\otimes a_{n+1}\mapsto a_0g(a_1\otimes\cdots\otimes a_n)a_{n+1}].\] We can assign $(Hom_k(A^{\otimes n},M))_{n\in\mathbb{N}}$ with differentials \[(\delta^n:Hom_k(A^{\otimes n},M)\rightarrow Hom_k(A^{\otimes n+1},M))_{n\in\mathbb{N}}\]
such that
\begin{align}
\begin{split}
    \delta^n f(a_1\otimes\cdots\otimes a_{n+1})&=a_1f(a_2\otimes\cdots\otimes a_{n+1})\\
    &+\sum_{i=1}^{n}(-1)^if(a_1\otimes\cdots\otimes a_ia_{i+1}\otimes \cdots\otimes a_{n+1})\\
    &+(-1)^{n+1}f(a_1\otimes\cdots\otimes a_n)a_{n+1}
\end{split}
\end{align}
for any $f\in Hom_k(A^{\otimes n}, M)$, as illustrated above.\\
In particular, when $M=A$, $\delta(id_A)\in Hom_k(A\otimes A, A)$ is the usual multiplication of $A$.\newline
\newline
Interpreting Hochschild cochains over $k$ instead of $A^e$ can greatly simplify our calculation. By virtue of $(4)$ we can analyze the behavior of Hochschild cochains in low degrees.
\subsection{Hochschild cochains in low degrees}
The differentials of Hochschild cochain complex $C^{\bullet}(A,M)$ in low degrees are given by
\[M\overset{\delta^0}{\longrightarrow}Hom_k(A,M)\overset{\delta^1}{\longrightarrow}Hom_k(A\otimes_kA,M)\overset{\delta^2}{\longrightarrow}Hom_k(A\otimes_kA\otimes_kA,M),\]
where
\begin{align}
    M\ni m&\mapsto [a\mapsto am-ma]\\
    Hom_k(A,M)\ni f&\mapsto [a\otimes b\mapsto af(b)-f(ab)+f(a)b]\\
    \begin{split}
    Hom_k(A\otimes A,M)\ni g&\mapsto [a\otimes b\otimes c\mapsto\\
    &ag(b\otimes c)-g(ab\otimes c)+ g(a\otimes bc)-g(a\otimes b)c]
    \end{split}
\end{align}
Note that $B_1$ in the definition of star product belongs to $C^2(A,A)$. The associativity of $B_1$ is equivalent to $\delta^2(B_1)=0$.\\
From $(5)$ we see that $HH^0(A,M)=\{m\in M|am=ma, \forall a\in A\}$. In particular, $HH^0(A)=Z(A)$, the center of $A$.\newline
\begin{definition}
Given a $k$-algebra $A$ and $M$ an $A$-bimodule. A $k$-linear map $f:A\rightarrow M$ is a \textit{$k$-derivation} if it satisfies Leibniz rule:
\begin{align*}
f(ab)=af(b)+f(a)b & & \forall a,b\in A.
\end{align*}
That is, $\delta^1(f)=0$.
The set of all $k$-derivations from $A$ to $M$ form a $k$-module, which is denoted by $Der(A,M)$; \\
Derivations of the form $a\mapsto am-ma$ are called \textit{inner derivations}, they form a sub $k$-module $InnDer(A,M)\subseteq Der(A,M)$;\\
The \textit{outer derivation} is the quotient space $Der(A,M)/InnDer(A,M)$, denoted $OutDer(A,M)$.
\end{definition}
It follows that $HH^1(A,M)=ker(\delta^1)/Im(\delta^0)$ is the $k$-module $OutDer(A,M)$. When $A$ is commutative, $InnDer(A)=0$. In that case,
\begin{align}
    HH^1(A)\cong Der(A).
\end{align}
\begin{remark}
$Der(A)$ naturally carries a Lie algebra structure: for $D_1, D_2\in Der(A)$, $[D_1,D_2]:=D_1\circ D_2-D_2\circ D_1$ is still a derivation. Being a commutator, $[-,-]$ is a Lie bracket.
\end{remark}
As an example, when $A=C^\infty(N)$ for some smooth manifold $N$, from a differential geometry course we know that $Der(A)$ consists of all ($C^\infty$)vector fields on $N$.
\subsection{Gerstenhaber bracket}
In practice, we are mainly concerned with Hochschild cochain $C^{\bullet}(A)$ for a $k$-algebra $A$ and its cohomology groups. So suppose $A$ is an associative algebra over $k$. We can regard $C^{\bullet}(A)$ as a graded $k$-
module with its natural grading
$C^{\bullet}(A)=\bigoplus\limits_{n=0}^{\infty} C^n(A)$.
\begin{definition}
     A \textit{graded Lie algebra} (GLA) is a $\mathbb{Z}$-graded vector space \[V=\bigoplus_{n=0}^{\infty} V_n,\] together with a bilinear map $[-,-]:V\times V\rightarrow V$ such that
     \begin{enumerate}
         \item \textit{respect grading:} $[V_i,V_j]\subseteq V_{i+j}$ ;
         \item \textit{graded skew-symmetry:} $[u,v]=-(-1)^{mn}[v,u]$, where $u$, $v$ are homogeneous elements in $V$ of degrees $m$, $n$, respectively;
         \item \textit{satisfies graded-Jacobi identity:}
         \begin{align*}
    (-1)^{mp}[u,[v,w]]+(-1)^{nm}[v,[w,u]]+(-1)^{pn}[w,[u,v]]=0,
         \end{align*}
where $u$, $v$, $w$ are homogeneous elements in $V$ of degrees $m$, $n$, $p$, respectively.
     \end{enumerate}
     $[-,-]$ is called a \textit{graded Lie bracket} on $V$.
\end{definition}
We can associate $C^{\bullet}(A)$ with a graded Lie bracket. Since $C^{\bullet}$ is generated by homogeneous element, we only need to define the bracket over homogeneous elements.
\begin{definition}
Let $f\in C^m(A,M)\cong Hom_k(A^{\otimes m},M)$, $g\in C^n(A,M)\cong Hom_k(A^{\otimes n},M)$, for $1\leq i\leq m+1$ define the \textit{$i^{th}$ circle product} by
\begin{align}
f\circ_i g:=f(id_A^{\otimes i-1}\otimes g\otimes id_A^{\otimes m-i})\in Hom_k(A^{\otimes m+n-1},M)\cong C^{m+n-1}(A).
\end{align}
Then define the circle product by
\begin{align}
    f\circ g:=\sum_{i=1}^m(-1)^{(i-1)(n+1)}f\circ_i g.
\end{align}
Finally, define the \textit{Gerstenhaber bracket} by
\begin{align}
    [f,g]:=f\circ g-(-1)^{(m-1)(n-1)}g\circ f
\end{align}
\end{definition}
\begin{figure}[ht]
    \centering
    \includegraphics[width=5cm]{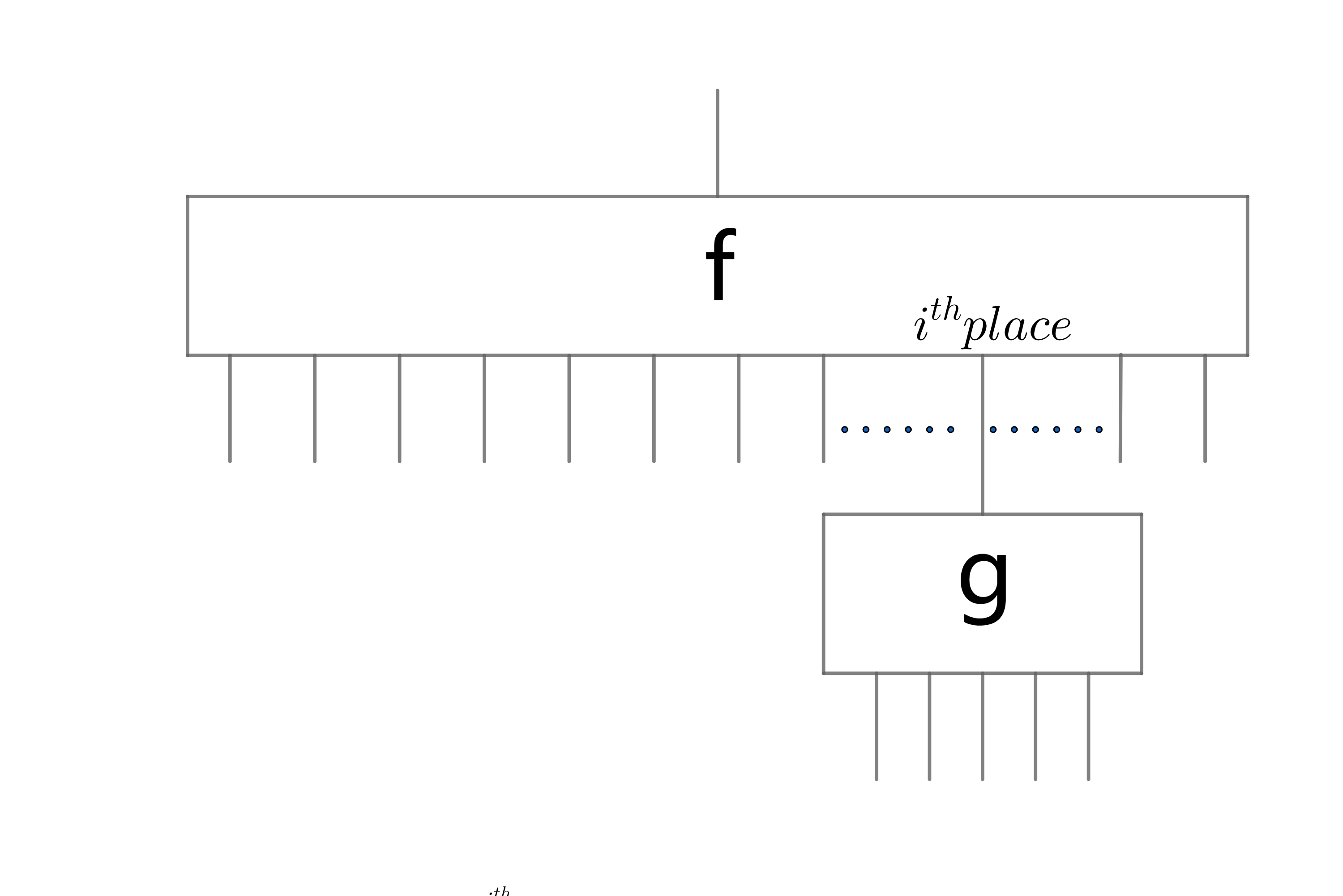}
    \caption{A graphic representation of $f\circ_i g$}
    \label{f circ g}
\end{figure}
\begin{remark}
The circle product $\circ$ defined above equips $C^{\bullet}(A)$ with \textit{pre-Lie structure}. Note that when $n=m=1$ the circle product $(10)$ is indeed composition of functions. However, in general, the circle product is not associative.
\end{remark}
The bracket $[-,-]$ in the above definition
gives a degree $0$ map from $C^{\bullet}(A)\times C^{\bullet}(A)$ to $C^{\bullet}(A)[-1]$. Applying a $+1$ shifting to both argument spaces
\begin{gather*}
    [-,-]:C^{\bullet}(A)\times C^{\bullet}(A)\longrightarrow C^{\bullet}(A)[-1]\\
               \Updownarrow\\
    [-,-]:C^{\bullet}(A)[1]\times C^{\bullet}(A)[1]\longrightarrow C^{\bullet}(A)[1]
\end{gather*}
gives a degree $0$ map from $C^{\bullet}(A)[1]\times C^{\bullet}(A)[1]$ to $C^{\bullet}(A)[1]$. So $[-,-]$ is a degree $0$ bracket \textit{up to a shift}.
\newline
$[-,-]$ is indeed graded skew-symmetric: for $f\in C^m(A)$, $g\in C^n(A)$, we have
\begin{align*}
-(-1)^{(m-1)(n-1)}[g,f]&=-(-1)^{(m-1)(n-1)}(g\circ f-(-1)^{(m-1)(n-1)}f\circ g)\\
&=f\circ g-(-1)^{(m-1)(n-1)}g\circ f=[f,g].
\end{align*}
After a shift by degree $1$ we have the graded skew-symmetry for a GLA.
One can verify $[-,-]$ satisfies graded-Jacobi identity by expanding the formula into $i^{th}$ circle products, see \cite{Gerstenhaber1963TheCS}.
Therefore,
\begin{lemma}
$(C^{\bullet}(A), [-,-])$ is a graded Lie algebra (up to a shift).
\end{lemma}
There is a link between Gerstenhaber bracket on $C^{\bullet}(A)$, the differential $\delta$, and the associative multiplication $\cdot$ in $A$. Let $\mu\in C^2(A)\cong Hom_k(A\otimes A, A)$ be the $2$-cochain corresponding to $\cdot$. Suppose $f\in C^n(A)$, the first and the third term in $RHS$ of $(4)$ gives $(-1)^{(2-1)(n-1)}(-\mu\circ f)$, while the second term gives $f\circ-\mu$. Thus
\begin{align}
    \delta(f)=[f,-\mu],
\end{align}
which then holds for any $f\in C^{\bullet}(A)$.\\
From this we can show that Gerstenhaber bracket is 'compatible' with the differential. For $f,g\in C^{\bullet}(A)$ where $g$ is homogeneous of degree $m$, we have
\begin{align*}
    \delta[f,g]&=[[f,g],-\mu]=(-1)^{m-1}[[f,-\mu],g]+[f,[g,-\mu]] & \textit{(by Jacobi identity)}\\
    &=(-1)^{m-1}[\delta(f),g]+[f,\delta g], \\
    &=(-1)^m[\delta(f),g]+[f,\delta g] & \textit{(after shifting)}.
\end{align*}
In defining Gerstenhaber bracket we regard $C^{\bullet}(A)$ as graded $k$-vector space only. With the same construction, for a $k$-vector space $V$, we can associate a Gerstenhaber bracket $[-,-]$ on $\bigoplus\limits_{n=0}^{\infty} Hom_k(V^{\otimes n},V)$.
\begin{definition}
$\mu\in Hom_k(V\otimes V, V)$ is an \textit{associative product on $V$} if for any $a,b,c\in V$,
\[\mu(a,\mu(b,c))=\mu(\mu(a,b),c).\]
\end{definition}
\begin{lemma}
$\mu\in Hom_k(V\otimes V, V)$ defines an associative product on $V$ if and only if $[\mu,\mu]=0$, where $[-,-]$ is the Gerstenhaber bracket.
\begin{proof}
By definition of Gerstenhaber bracket with $m=n=2$,
\begin{align*}
    \frac{1}{2}[\mu,\mu]&=\frac{1}{2}(\mu\circ\mu-(-1)^3\mu\circ\mu)=\mu\circ\mu\\
    &=\mu\circ_1\mu -\mu\circ_2\mu\\
    &=\mu(\mu\otimes id_V)-\mu(id_V\otimes\mu)=0 & \iff & \text{$\mu$ is associative}.
\end{align*}
\end{proof}
\end{lemma}
\begin{theorem}
Let $A$ be an associative algebra with underlying vector space $V$ and multiplication $\mu$. Let $\nu\in C^2(A)\cong Hom_k(A\otimes A, A)$ be a Hochschild $2$-cochain. Then $\mu+\nu\in C^2(A)$ is associative if and only if \[\delta(\nu)+\frac{1}{2}[\nu,\nu]=0\]
\begin{proof}
By above lemma, $\mu+\nu$ is associative if and only if
\begin{align*}
    0&=\frac{1}{2}[\mu+\nu,\mu+\nu]\\
    &=\frac{1}{2}([\mu,\mu]+[\nu,\nu]+[\mu,\nu]+[\nu,\mu])\\
    &=\frac{1}{2}[\nu,\nu]+\delta\nu
\end{align*}
\end{proof}
\end{theorem}

We now have a process to identify deformed products of an associative algebra $(A,\mu)$ on $A[[\hbar]]$:
\begin{enumerate}
    \item Extend the product $\mu$ to $A[[\hbar]]$ by defining\footnote{or $(\sum\limits_{i\geq 0}a_i\hbar^i,\sum_{j\geq 0}b_j\hbar^j)\mapsto \sum\limits_{r\geq 0}\sum\limits_{0\leq l\leq r}(a_{r-l}b_l)\hbar^r.$}
    \[(\sum_{i\geq 0}a_i\hbar^i,\sum_{j\geq 0}b_j\hbar^j)\mapsto a_0b_0,\]
    and call this bilinear map $\mu_0$, note that $\mu_0$ is an associative product on $A[[\hbar]]$.
    \item By definition, deformed products of $\mu$ on $A[[\hbar]]$ are of the form \[\mu+\sum\limits_{k\geq 1}\hbar^k\mu_k,\]
    where each $\mu_k$ is in $C^2(A)$.
    \item The theorem says that such deformed products $\mu+\sum\limits_{k\geq 1}\hbar^k\mu_k$ of $\mu$ are precisely given by the relation
    \[d(\sum_{k\geq 1}\hbar^k\mu_k)+\frac{1}{2}[\sum_{k\geq 1}\hbar^k\mu_k,\sum_{k\geq 1}\hbar^k\mu_k]=0.\]
\end{enumerate}

\newpage
\section{DGLAs and Maurer-Cartan elements}
\subsection{Differential graded Lie algebras}
\begin{definition}
A \textit{differential graded Lie algebra} (\textit{DGLA}) is a graded Lie algebra $L$ together with a degree $1$ linear map $d:L\rightarrow L$ which is
\begin{enumerate}
    \item a differential, $i.e.$ $d^2=0$;
    \item a degree $1$ derivation, $i.e.$ $d([f,g])=(-1)^{m}[d(f),g]+[f,d(g)]$ for $f,g\in L$, $g$ homogeneous of degree $m$.
\end{enumerate}
\end{definition}
\textbf{Examples}
\begin{enumerate}
    \item Up to a shifting, Hochschild cochain $(C^{\bullet}(A), \delta)$ for an associative algebra $A$ together with Gerstenhaber bracket $[-,-]$ is a DGLA.
    \item Suppose $M$ is a manifold and $A=C^{\infty}(M)$. For each $n\in \mathbb{N}$ we take $D^n(A)\subset C^n(A)$ to be the ($k$-)submodule of $C^n(A)$ consisting of linear polydifferential operators. By definition  $\bigoplus\limits_{n=0}^{\infty} D^n(A)$ is closed under both $\delta$ and $[-,-]$. Hence $(\bigoplus\limits_{n=0}^{\infty} D^n(A), \delta, [-,-])$ is a subchain of $C^n(A)$ and is also a DGLA. We call it $D_{poly}$.
    \item If $(L=\bigoplus\limits_{n=0}^{\infty}L^n,[-,-],d)$ is a DGLA (over $k$), then it induces a DGLA structure on $L[[\hbar]]:=L\otimes_kk[[\hbar]]$ over $k[[\hbar]]$. The $n^{th}$ component of $L[[\hbar]]$ is $L^n\otimes_kk[[\hbar]]$, with differential $d_{\hbar}:=d\otimes_k id_{k[[\hbar]]}$ and the graded-Lie bracket is defined by
    \begin{align*}
        [\sum_{i\geq 0}f_i\hbar^i,\sum_{j\geq 0}g_j\hbar^j]:=\sum_{k\geq 0}\hbar^k\sum_{i+j=k}[f_i,g_j].
    \end{align*}
\end{enumerate}
In the sequel we are going to shift the DGLA $C^{\bullet}(A)$ (resp. $D_{poly}$) so that $C^2(A)$ (resp. $D^2(A)$) lives in degree $1$ position. We still denote the resulting DGLA as $C^{\bullet}(A)$ (resp. $D_{poly}$).\footnote{Note we have shifted the chain "by a place to the right" so $D^n_{poly}= D^{n+1}(A)$.} \\
Another important example of DGLA is the chain of \textit{polyvector fields}. Recall that $Der(A)$ is equipped with a Lie bracket $[-,-]$ defined as commutator. We can extend this bracket to $\bigoplus\limits_{n=0}^{\infty}\bigwedge^{n+1}Der(A)$:
\begin{definition}
Let $\eta_0\wedge\cdots\wedge\eta_n\in\bigwedge^{n+1}Der(A)$ and $\xi_0\wedge\cdots\wedge\xi_m\in\bigwedge^{m+1}Der(A)$ be polyvector fields. The \textit{Schouten-Nijenhuis bracket} is given by
\begin{equation*}
[\eta_0\wedge\cdots\wedge\eta_n,\xi_0\wedge\cdots\wedge\xi_m]_{SN}:=\sum_{\substack{0\leq i\leq n\\0\leq j\leq m}}(-1)^{i+j}[\eta_i,\xi_j]\wedge\eta_0\cdots\wedge\hat{\eta_i}\wedge\cdots\wedge\eta_n\wedge\xi_0\wedge\cdots\wedge\hat{\xi_j}\wedge\cdots\wedge\xi_m,
\end{equation*}
where $\hat{}$ stands for omission.
\end{definition}
\begin{remark}
Taking $m=n=0$, $[-,-]_{SN}$ coincide with the usual Lie bracket on $Der(A)$. In fact, this bracket is the \textit{universal extension} of the Lie bracket on $Der(A)$ to a graded version.
\end{remark}
For $A=C^{\infty}(\mathbb{R}^n)$ for some $n\in\mathbb{N}^* $, there is another way to construct a graded-Lie beacket on  $\bigoplus\limits_{n=0}^{\infty}\bigwedge^{n+1}Der(A)$:\newline
In local coordinate $(x_1,\cdots,x_n)$, a homogeneous element has the form
\[\Pi=\sum_{i_1<\cdots<i_p}\Pi^{i_1,\cdots,i_p}\frac{\partial}{\partial x_{i_1}}\wedge\cdots\wedge\frac{\partial}{\partial x_{i_p}}\]
Write every $\frac{\partial}{\partial x_i}$ as $\zeta_i$ and regard them as formal and \textit{odd} variables (odd in the sense that $\zeta_i\zeta_j:=\frac{\partial}{\partial x_i}\wedge\frac{\partial}{\partial x_j}$ so $\zeta_i\zeta_j=-\zeta_j\zeta_i$). $\Pi$ is then a function in $(x_1,\cdots,x_n,\zeta_1,\cdots,\zeta_n)$. We define, for $\Pi_1\in\bigwedge^{p_1}Der(A),\Pi_2\in\bigwedge^{p_2}Der(A)$,
\[\Pi_1\bullet\Pi_2:=\sum_{i=1}^{n}\frac{\partial \Pi_1}{\partial \zeta_i}\frac{\partial\Pi_2}{\partial x_i}\]
and
\[[\Pi_1,\Pi_2]:=\Pi_1\bullet\Pi_2-(-1)^{(p_1-1)(p_2-1)}\Pi_2\bullet\Pi_1\]
It is indeed an extension of the standard Lie bracket on $Der(A)=\mathfrak{X}(M)$:\\
For $a=\sum\limits_{i=1}^na_i\frac{\partial}{\partial x_i}=\sum\limits_{i=1}^na_i\zeta_i$ and $b=\sum\limits_{i=1}^nb_i\frac{\partial}{\partial x_i}=\sum\limits_{i=1}^nb_i\zeta_i$ then
\begin{align*}
    [a,b]&=\sum_{i=1}^na_i(\sum_{j=1}^n\frac{\partial b_j}{\partial x_i}\frac{\partial}{\partial x_j})
    -\sum_{i=1}^nb_i(\sum_{j=1}^n\frac{\partial a_j}{\partial x_i}\frac{\partial}{\partial x_j})\\
    &=\sum_{i=1}^n\frac{\partial a}{\partial\zeta_i}\frac{\partial b}{\partial x_i}
    -\sum_{i=1}^n\frac{\partial b}{\partial\zeta_i}\frac{\partial a}{\partial x_i}
\end{align*}
In fact, this Lie bracket on $\bigoplus\limits_{n=0}^{\infty}\bigwedge^{n+1}Der(A)$ is identical to the Schouten-Nijenhuis bracket $[-,-]_{SN}$. The graded skew-symmetry follows, because it is a commutator of $\bullet$. For graded-Jacobi identity and more generalities about $[-,-]_{SN}$, see \cite{article} and \cite{2005}.
\newline
Equipping $\bigoplus\limits_{n=0}^{\infty}\bigwedge^{n+1}Der(A)$ with trivial differential $d=0$ turns it into a cochain complex which is a DGLA:
\begin{definition}
$T_{poly}$ is the DGLA $\bigwedge^{\bullet}Der(A)$ equipped with differential $d=0$ and the Schouten-Nijenhuis bracket. It is shifted so that $\wedge^2Der(A)$ is in the degree 1 position.
\end{definition}
Indeed, since $d=0$ the compatibility with $d$ is trivial.\newline\newline
Suppose $M$ is a $C^{\infty}$ manifold and $A=C^{\infty}(M)$. Then the degree 1 component of the corresponding $T_{poly}$ are bivector fields on $M$. Each bivector field gives a bracket $\{-,-\}:A\times A\rightarrow A$ which is skew-symmetric and satisfies Leibniz rule.
\begin{theorem}
Let $\Pi\in\wedge^2Der(A)$ be a bivector field on a manifold $M$ of dimension $n$. It defines a Poisson bracket on $A$ if and only if \[d\Pi+\frac{1}{2}[\Pi,\Pi]_{SN}=0.\]
\begin{proof}
Since $d=0$, it is equivalent to show that
\begin{center} $\Pi$ is Poisson $\iff$ $[\Pi,\Pi]_{SN}=0$.\end{center}
In local coordinates, we can write $\Pi=\sum\limits_{1\leq i,j\leq n}\Pi^{ij}\frac{\partial}{\partial x_i}\wedge\frac{\partial}{\partial x_j}$. Direct calculation gives
\begin{align*}
    [\Pi,\Pi]_{SN}&=[\sum_{1\leq i,j\leq n}\Pi^{ij}\frac{\partial}{\partial x_i}\wedge\frac{\partial}{\partial x_j},\sum_{1\leq k, l\leq n}\Pi^{kl}\frac{\partial}{\partial x_k}\wedge\frac{\partial}{\partial x_l}]_{SN}\\
    &=\sum_{1\leq i,j,k,l\leq n}\Pi^{ij}\frac{\partial \Pi^{kl}}{\partial x_j}\frac{\partial}{\partial x_i}\wedge\frac{\partial}{\partial x_k}\wedge\frac{\partial}{\partial x_l}=0
\end{align*}
$\iff$ for all $f,g,h\in C^{\infty}(M)$,
\begin{align*}
    0=\sum_{1\leq i,j,k,l\leq n}\Pi^{ij}\frac{\partial \Pi^{kl}}{\partial x_j}\frac{\partial}{\partial x_i}\wedge\frac{\partial}{\partial x_k}\wedge\frac{\partial}{\partial x_l}(df,dg,dh).
\end{align*}
Expanding the $RHS$ gives a (nonzero) multiple of
\begin{align*}
    &\sum_{1\leq i,j,k,l\leq n}
    ((\Pi^{ij}\frac{\partial \Pi^{kl}}{\partial x_j})
    (\frac{\partial f}{\partial x_i}\frac{\partial g} {\partial x_k}\frac{\partial h}{\partial x_l}
    +\frac{\partial g}{\partial x_i}\frac{\partial h}{\partial x_k}\frac{\partial f}{\partial x_l}
    +\frac{\partial h}{\partial x_i}\frac{\partial f}{\partial x_k}\frac{\partial g}{\partial x_l}))\\
    &=\{\{f,g\},h\}+\{\{g,h\},f\}+\{\{h,f\},g\}.
\end{align*}
\end{proof}
\end{theorem}
\subsection{Maurer-Cartan elements}
\begin{definition}
Suppose $(L,[-,-],d)$ is a DGLA. An element $s$ of degree $1$ is said to be \textit{Maurer-Cartan} if it satisfies the \textit{Maurer-Cartan equation}:
\[ds+\frac{1}{2}[s,s]=0\]
We denote $MC(L)$ for the set of Maurer-Cartan elements in $L$.
\end{definition}
Let us look at Maurer-Cartan elements of some variations of DGLAs.
\begin{itemize}
    \item By Theorem 22, for a manifold $M$, $MC(T_{poly})=MC(\bigwedge^{\bullet}Der(C^{\infty}(M)))$ corresponds to the set of all Poisson brackets on $M$.
    \item If we let $A:=C^{\infty}(M)$, then $MC(\hbar C^{\bullet}[[\hbar]])$ corresponds to all products on $A[[\hbar]]$ deforming the standard pointwise product $\mu$ on $A$. In the language of deformation quantization, we want the star product given by bilinear bidifferential operators, which means that we are looking for the set \[\{\sum_{k\geq 1}\hbar^k\mu_k\in MC(\hbar C^{\bullet}[[\hbar]])\text{ }|\text{ }\mu_k\in D^2(A)\},\]but this is precisely $MC(\hbar D_{poly}[[\hbar]])$.
    \item A combined argument of the above two situations shows that the formal Poisson structures on $M$ are given by $MC(\hbar T_{poly}[[\hbar]])$.
    \item Finally, suppose $L$ is a DGLA and $\tilde{L}:=L\otimes (\hbar)=\hbar L[[\hbar]]$(see Example 3 in $\S4.1$). The degree $n$ elements in $\tilde{L}$ are of form $\sum\limits_{i\geq 1}f_i\hbar^i$ where each $f_i\in L^n$. Such element is Maurer-Cartan if and only if
    \begin{align}
       df_k+\frac{1}{2}\sum_{i+j=k}[f_i,f_j]=0\tag{$MC_k$}
    \end{align}
    for each $k\geq 1$, by looking at the coefficient of $\hbar^k$.
\end{itemize}
Since the (intended) identification
\begin{center}
    $\{$star products$\}\longleftrightarrow\{$(formal) Poisson structures$\}$
\end{center}
are up to gauge equivalence on both sides, our next task is to construct the gauge group action on Maurer-Cartan elements.
\subsection{Gauge group actions}
We assume $L$ is a DGLA whose differential $d$ has the form $d=[-,\theta]$ for some $\theta\in L^1$ satisfying $[\theta,\theta]=0$. Since the DGLAs $T_{poly}$, $C^{\bullet}(A)$ (see $(12)$) and $D_{poly}$ all have such $\theta$ in their first degree components, it suffices for our needs. \\
Let $\tilde{L}:=L\otimes (\hbar)=\hbar L[[\hbar]]$. Then the differential $\tilde{d}$ in $\tilde{L}$ is also of the form $\tilde{d}=[-,\theta]$, although $\theta\notin\tilde{L}$. \\
Furthermore, We can regard $\tilde{L}^0=L^0\otimes (\hbar)$ as a Lie-subalgebra of $\tilde{L}$, and $\tilde{L}^1=L^1\otimes (\hbar)$ as a (Lie-)representation of $\tilde{L}^0$, with the action of $\alpha\in \tilde{L}^0$ on $l\in\tilde{L}^1$ given by:
\begin{align}
    \alpha\cdot l&:=l+[\alpha,\theta+l]+\frac{1}{2!}[\alpha,[\alpha,\theta+l]]+\frac{1}{3!}[\alpha,[\alpha,[\alpha,\theta+l]]]+\cdots\\
    &=exp(\alpha)(\theta+l)exp(-\alpha)-\theta
\end{align}
Under our assumption, $(13)$ can be expanded in the form
\begin{align}
    \alpha\cdot l=l+d\alpha+[\alpha,l]+\textit{higher order brackets}
\end{align}
Note that since $[-,-]$ is a degree $0$ bracket and each $\alpha$ is a multiple of $\hbar$, each component of $(13)$ lives in $\tilde{L}^1$ so the action is well-defined.
\begin{definition}
For DGLA $\tilde{L}$ as mentioned above, let $\tilde{\Gamma}:=exp(\tilde{L}^0)$ be the gauge group. The group action of $\Gamma$ on $\tilde{L}^0$ is given by
\[\chi\cdot l:=\chi (l+\theta)\chi^{-1}-\theta\] for $\chi\in\tilde{\Gamma}$ and $l\in\tilde{L}^1$.
\end{definition}
\begin{lemma}
The action of $\tilde{\Gamma}$ on $\tilde{L}^0$ preserves $MC(\tilde{L})$.
\begin{proof}
The condition of Maurer-Cartan on $l$ is
\begin{align*}
    0&=2\tilde{d}l+[l,l]\\
    &=2[l,\theta]+[l,l]+[\theta,\theta] & \text{($[\theta,\theta]=0$)}\\
    &=[l,\theta]+[\theta,l]+[\theta,\theta]+[l,l] \\
    &=[\theta+l,\theta+l],
\end{align*}
where the third equality comes from graded skew-symmetric of $[-,-]$ and the degrees of $\theta$ and $l$.\\
Similarly, the Maurer-Cartan equation for $\alpha\cdot l$ is
\begin{align*}
[exp(\alpha)(\theta+l)exp(-\alpha),exp(\alpha)(\theta+l)exp(-\alpha)]=0,
\end{align*}
which can be rearranged as
\[exp(\alpha)[\theta+l,\theta+l]exp(-\alpha)=0\]
\end{proof}
\end{lemma}
Thanks to the lemma, for a DGLA $L$ we can define
\[Def(L):=MC(\hbar L[[\hbar]])/\tilde{\Gamma}.\]
When $L=T_{poly}$ over a Poisson manifold $(M,\Pi)$, $Def(L)$ represents equivalence classes of Formal Poisson structures over $M$ modulo gauge group action, as seen in $\S2.3$. \\
For $L=C^{\bullet}(A)$, $Def(L)$ gives equivalence classes of star products:
\begin{lemma}
Given an associative algebra $(A,\mu)$, two star products \[-*_1-:=\mu+\sum\limits_{i\geq 1}\hbar^i\mu_i\] and \[-*_2-:=\mu+\sum\limits_{i\geq 1}\hbar^i\mu'_i\] are gauge equivalent if and only if their associated Maurer-Cartan elements $\sum\limits_{i\geq 1}\hbar^i\mu_i$ and $\sum\limits_{i\geq 1}\hbar^i\mu'_i$ are gauge equivalent.
\begin{proof}
(Compare $\S2.3$) In this case $L^0=Hom_k(A,A)$, an element $\alpha$ of $\tilde{L}^0$ has the form
\[\alpha=f_1\hbar+f_2\hbar^2+\cdots, \hfill{f_i\in Hom_k(A,A)}.\]
Then
\[exp(\alpha)=id_A+\alpha+\frac{1}{2!}\alpha^2+\cdots:A[[\hbar]]\rightarrow A[[\hbar]]\] is an invertible $k[[\hbar]]$-linear map. The gauge equivalence of $\sum\limits_{i\geq 1}\hbar^i\mu_i$ and $\sum\limits_{i\geq 1}\hbar^i\mu'_i$ means, there exists $\chi=exp(\alpha)\in exp(\tilde{L}^0)=\tilde{\Gamma}$ such that, as functions on $A[[\hbar]]\otimes A[[\hbar]]$,
\begin{align*}
    \sum\limits_{i\geq 1}\hbar^i\mu'_i&=\chi\cdot\sum\limits_{i\geq 1}\hbar^i\mu_i=\chi(\mu+\sum\limits_{i\geq 1}\hbar^i\mu_i)\chi^{-1}-\mu,
\end{align*}
so
\[ \mu+\sum\limits_{i\geq 1}\hbar^i\mu'_i=\chi(\mu+\sum\limits_{i\geq 1}\hbar^i\mu_i)\chi^{-1},\]
which means
\[ (\mu+\sum\limits_{i\geq 1}\hbar^i\mu'_i)(exp(\alpha)\otimes exp(\alpha))=exp(\alpha)(\mu+\sum\limits_{i\geq 1}\hbar^i\mu_i)
\]
as functions from $A[[\hbar]]\otimes A[[\hbar]]$ to $A[[\hbar]]$. Hence $-*_1-$ and $-*_2-$ are gauge equivalent via $D=exp(\alpha)$.\\
Conversely, every invertible $D:A[[\hbar]]\rightarrow A[[\hbar]]$ such that \[D=id+D_1\hbar+D_2\hbar^2+\cdots\]
is of the form $exp(\alpha)$ for some $\alpha\in Hom_k(A,A)\otimes (\hbar)$\footnote{Because $\hbar$ is nilpotent.}, as learned from a Lie group course. Moving backward shows that the gauge equivalence of $-*_1-$ and $-*_2-$ implies gauge equivalence of
$\sum\limits_{i\geq 1}\hbar^i\mu_i$ and $\sum\limits_{i\geq 1}\hbar^i\mu'_i$.
\end{proof}
\end{lemma}
For a Poisson manifold $(M,\Pi)$ with $A:=C^{\infty}(M)$, it follows immediately from Lemma 26 that $Def(D_{poly})$ gives equivalence classes of deformation quantizations.\\
To prove the main theorem, we want to identify $Def(T_{poly})$ and $Def(D_{poly})$. The ideal strategy would be as follows: first, to find an isomorphism between DGLAs $T_{poly}$ and $D_{poly}$ which preserves the Lie bracket. Second, to extend this to an isomorphism between $T_{poly}\otimes(\hbar)$ and $D_{poly}\otimes(\hbar)$. Third, to use this extended isomorphism to identify $MC(T_{poly}\otimes(\hbar))$ with $MC(D_{poly}\otimes(\hbar))$ modulo gauge group actions.\\
For this moment if we regard $T_{poly}$ and $D_{poly}$ as chain complexes of $k$-($\mathbb{R}$-)modules, on the $0^{th}$ components we have a $k$-($\mathbb{R}$-)linear map, naturally defined as
\[\Psi_0:T^0_{poly}=Der(A)=\mathfrak{X}(M)\rightarrow D^1(A):\xi\mapsto [f\mapsto \xi f].\]
This (naturally) generalizes to\footnote{Regarded as wedge products of covectors $\xi_i$'s on the vector space $A$, for example.}, for $n\in\mathbb{N}$,
\begin{align*}
\Psi_n:\xi_0\wedge\cdots\wedge\xi_n\mapsto[(f_0\otimes\cdots\otimes f_n)\mapsto\frac{1}{(n+1)!}\sum_{\sigma\in S_{n+1}}sgn(\sigma)\prod_{i=0}^{n}\xi_{\sigma(i)}(f_i)],
\end{align*}
where $sgn(\sigma)$ is the sign of the $n+1$-permutation $\sigma\in S_{n+1}$.\\
One can check that $\Psi=(\Psi_i)_{i\in\mathbb{N}}:T_{poly}\rightarrow D_{poly}$ is a map of chains using $(4)$ and the fact that the differential $d$ in $T_{poly}$ is $0$.  \\
\begin{remark}
Unfortunately, $\Psi$ fails to preserve the Lie bracket, as one can calculate, for example,
\[\Psi([\xi_1\wedge\xi_2,\xi_3\wedge\xi_4]_{SN})(f\otimes g\otimes h)=\frac{1}{6}(\textit{sum of elements of form }\pm\xi_i(\xi_j(f)\xi_k(g)\xi_l(h))),\]
while
\[[\xi_1\wedge\xi_2,\xi_3\wedge\xi_4](f\otimes g\otimes h)=\frac{1}{4}(\textit{sum of elements of form }\pm\xi_i(f)\xi_j(\xi_k(g)\xi_l(h))).\]
\end{remark}
It turns out that, though $\Psi$ does not preserve Lie bracket, it is a quasi-isomorphism ($i.e.$ it induces isomorphism between cohomology groups of the chains) between $T_{poly}$ and $D_{poly}$.
\begin{theorem}
(\textit{Hochschild-Kostant-Rosenberg}, HKR \cite{Hochschild1962DifferentialFO})
For a smooth manifold $M$ with $A:=C^{\infty}(M)$, its associated chain complexes $T_{poly}=\bigwedge^{\bullet}(Der(A))$ and $D_{poly}\subset C^{\bullet}(A)$ are quasi-isomorphic via $\Psi$.
\end{theorem}
In appendix B we give a proof of HKR in the version of
\[\Phi: \bigwedge\nolimits^{\bullet}(Der(A))\rightarrow C^{\bullet}(A),\]
when $A$ is a \textit{smooth} commutative noetherian $k$-algebra and $\Phi$ is given exactly the same as $\Psi$ above. The version of quasi-isomorphism between $T_{poly}$ and $D_{poly}$ over a smooth manifold can be found in \cite{Kontsevich_2003}.\\
Note that in the case of $A=C^{\infty}(M)$, the isomorphism induced by $\Phi$ on the level of cohomology groups between $\bigwedge^{\bullet}(Der(A))$ and $C^{\bullet}(A)$ does not directly imply the isomorphism on the level of cohomology groups between
$\bigwedge^{\bullet}(Der(A))=T_{poly}$ and $D_{poly}$, even $\Phi$ factors through $\Psi$.
\newpage
\section{$L_{\infty}$-algebras}
From the previous section (and appendix B) we have a quasi-isomorphism between $T_{poly}$ and $D_{poly}$ which fails to preserve the Lie brackets, not to mention their Maurer-Cartan elements. To resolve this problem, we should work in the context of $L_{\infty}$-algebras and $L_{\infty}$-morphisms. For a rigorous discussion one should begin with the notion of \textit{coassociative coalgebra} defined on certain graded vector spaces. However, since the general procedure of unwrapping a coalgebra structure is not particularly enlightening at this stage, we only summarize some key features that serve to build a connection for what we have done so far. Please find \cite{inbook} and \cite{Lada1994StronglyHL} for detailed construction of $L_{\infty}$-algebras.
\newline\newline
Throughout, for any homogeneous element $v$ in a graded vector space $V$, denote $|v|$ for the degree of $v$.
\begin{definition}
Let $V$ be a $\mathbb{Z}$-graded vector space. \\
The \textit{tensor algebra} over $V$ is the algebra
\[T(V):=\bigoplus_{n=0}^{\infty}V^{\otimes n}\]
with multiplication $\otimes$ and graded in the way that $v_1\otimes\cdots\otimes v_n$ has degree $\sum\limits^{n}_{i=1}m_i$, where  $v_1,\cdots, v_n\in V$ are homogeneous of degrees $m_1,\cdots,m_n$ respectively. We set $V^{\otimes 0}:=0$.  \\
The \textit{free-graded commutative associative algebra} over $V$ is the algebra
\[\bigwedge V:=T(V)/(x\otimes y-(-1)^{pq}y\otimes x),\]
where $x,y$ are homogeneous elements in $V$ of degrees $p,q$ respectively.
\end{definition}
Note that the multiplication and grading in $\bigwedge V$ are induced from the multiplication $\otimes$ and grading in $T(V)$.\\
We define, for each $n\in\mathbb{N}^*$ and $n$-permutation $\sigma\in S_n$, the \textit{Koszul sign} of $\sigma$, $\tau(\sigma)\in\{\pm 1\}$ satisfying
\[v_1\wedge\cdots\wedge v_n=\tau(\sigma)v_{\sigma(1)}\wedge\cdots\wedge v_{\sigma(n)}\]
for each $\sigma\in S_n$ and $v_1,\cdots,v_n\in\bigwedge V$,\\
and define the \textit{(graded) skew-symmetric Koszul sign} $\kappa$ as
\[\kappa(\sigma):=sgn(\sigma)\tau(\sigma).\]
$\kappa$ deserves its name, for suppose $\sigma\in S_n$ is the pemutation swapping $i$ and $j$, then $\tau(\sigma)=(-1)^{|v_i||v_j|}$ and $\kappa(\sigma)=-(-1)^{|v_i||v_j|}$.\\
We call a permutation $\sigma\in S_n$ an $(i,n-i)$-\textit{unshuffle} if $\sigma(1)<\cdots <\sigma(i)$ and $\sigma(i+1)<\cdots<\sigma(n)$. The set of all $(i,n-i)$-unshuffles is denoted as $S_{(i,n-i)}$. $S_{(i,n-i)}$ clearly forms a subgroup of $S_n$.  \\
Similarly, an $(i_1,\cdots,i_m)$-\textit{unshuffle} is a permutation $\sigma\in S_n$ where $n=i_1+\cdots +i_m$ such that the order is preserved within each block of length $i_1,\cdots,i_m$. By $S_{i_1+\cdots+i_m=n}$ we mean the union of subgroups in $S_n$ of the form $S_{i_1,\cdots,i_m}$ where $i_1+\cdots+i_m=n$.
\begin{definition}
An $L_{\infty}$-structure is a graded vector space $V$ together with a system of linear maps $(l_n:\otimes^n V\rightarrow V)_{n\in \mathbb{N}}$, each $l_n$ of degree $2-n$,
satisfying the following properties:
\begin{enumerate}
    \item \textit{(graded) anti-symmetry}: for each $n\in\mathbb{N}$, every permutation $\sigma\in S_n$ and homogeneous elements $v_1,\cdots,v_n\in V$,
    \begin{align*}
    l_n(v_{\sigma(1)},\cdots,v_{\sigma(n)})=\kappa(\sigma)l_n(v_1,\cdots,v_n);
    \end{align*}
    \item \textit{$L_n$}: For each $n\geq 1$ and homogeneous elements $v_1,\cdots,v_n\in V$,
    \begin{align*}
    \sum_{i+j=n+1}\sum_{\sigma\in S_{(i,n-i)}}\kappa(\sigma)l_j(l_i(v_{\sigma(1)},\cdots,v_{\sigma(i)}),v_{\sigma(i+1)},\cdots,v_{\sigma(n)})=0.
    \end{align*}
\end{enumerate}
\end{definition}
\begin{remark}
By linearity of each $l_n$, if $L=(V,(l_n)_{n\in\mathbb{N}})$, $L'=(V',(l'_n)_{n\in\mathbb{N}})$ are $L_{\infty}$-algebras, then $L\oplus L'$ is also an $L_{\infty}$-algebra. \\
If we set $l_k:=0$ for each $k>m$, every $L_m$ structure can be regarded as an $L_{\infty}$-algebra.
\end{remark}
We can investigate $L_n$ for small n's:
\begin{enumerate}
    \item $n=1$: In this case $L_1$ can be written as $l_1(l_1(v))=0$ for every $v\in V$, $i.e.$ $l_1$ is a degree $+1$ differential;
    \item $n=2$: The anti-symmetry condition implies that
    \[l_2(v_1,v_2)=-(-1)^{|v_1||v_2|}l_2(v_2,v_1)\]
     is a degree $0$-linear map which is graded skew-symmetric, whereas $L_2$ gives \[l_1(l_2(v_1,v_2))=l_2(l_1(v_1),v_2)+(-1)^{|v_1||v_2|}l_2(v_1,l_1(v_2)).\]
     If we write $[-,-]$ for $l_2(-,-)$ and $d(-)$ for $l_1(-)$, then $L_2$ becomes
     \[d[v_1,v_2]=[dv_1,v_2]+(-1)^{|v_1|}[v_1,v_2],\]
     which says that a DGLA is a special case of $L_{\infty}$-algebra with all $l_n=0$ for $n\geq 3$.
\end{enumerate}
The morphism between two $L_{\infty}$-algebras $L$ and $L'$ is given by unwrapping the definition of a morphism between \textit{differential graded coalgebras} (see \cite{Kajiura_2006} ):
\begin{definition}
A morphism of $L_{\infty}$-algebras from $L=(V,(l_n)_{n\in\mathbb{N}})$ to
$L'=(V',(l'_n)_{n\in\mathbb{N}})$ is a collection of graded skew-symmetric
multilinear maps
\[f:=(f_n:\otimes^n V\rightarrow V')_{n\geq 1},\]
with each $f_n$ is of degree $1-n$, satisfying
\begin{align*}
    &\sum_{\sigma\in S_{i_1+i_2=n}}(-1)^{\kappa(\sigma)}f_{1+i_2}(l_{i_1}(v_{\sigma(1)}\otimes\cdots\otimes v_{\sigma(i_1)})\otimes v_{\sigma(i+1)}\otimes \cdots\otimes v_{\sigma(n)})\\
    &=\sum_{\sigma\in S_{i_1+\cdots +i_m=n}}\frac{(-1)^{\kappa(\sigma)}}{m!}l_m'(f_{i_1}(v_{\sigma(1)}\otimes\cdots\otimes v_{\sigma(i_1)})\otimes\cdots\otimes f_{i_m}(v_{\sigma(n-i_m+1)}\otimes\cdots\otimes v_{\sigma(v_n)}))
\end{align*}
A morphism of $DGLAs$ is a morphism of $L_{\infty}$-algebras when both DGLAs are regarded as $L_{\infty}$-algebras.
\end{definition}
\begin{remark}
If $L'$ is a DGLA, the formula can be greatly simplified since $m=1\text{ or }2$ in this case (see \cite{Lada1994StronglyHL}).
\end{remark}
Given $h_1,h_2,\cdots,h_n\in\mathbb{N}$ and an $L_{\infty}$ morphism $f=(f_1,f_2,\cdots)$ in the above definition,
the restriction of $f_n$ on $V^{h_1}\otimes\cdots\otimes V^{h_n}\subset V^{\otimes n}$ is a linear map
\[f_{(h_1,\cdots,h_n)}:V^{h_1}\otimes\cdots\otimes V^{h_n}\rightarrow V^{h_1+\cdots+h_n+(1-n)}\]
satisfying the graded skew-symmetry property (from antisymmetry of $L_{\infty}$ structure on $V$):
\[f_{(h_1,\cdots,h_n)}(v_1\otimes\cdots\otimes v_n)=-(-1)^{h_ih_{i+1}}f_{(h_1,\cdots,h_{i+1}h_i,\cdots,h_n)}(v_1\otimes\cdots\otimes v_{i+1}\otimes v_i \otimes\cdots\otimes v_n)\]
where $v_j\in V^{h_j}$, $j=1,2,\cdots,n$. By abuse of notation we will write
\[f_n(v_1\wedge\cdots\wedge v_n):=f_{(h_1,\cdots,h_n)}(v_1\otimes\cdots\otimes v_n).\]
One can investigate the behavior of components $f_n$ of an $L_{\infty}$-morphism for small $n$:
\begin{enumerate}
    \item $n=1$: the condition can be written as \[f_1(l_1(v))=l_1'(f_1(v)),\] indicating that $f$ is compatible with the differentials;
    \item $n=2$: the condition can be written as \[f_1(l_2(v_1\otimes v_2))=l_2(f_1(v_1)\otimes f_1(v_2)),\] indicating that $f$ is compatible with the Lie brackets.
\end{enumerate}
Hence the $L_{\infty}$ morphisms are the kind of morphisms we want to work with. Given a morphism of DGLAs $f:L\rightarrow L'$
Points $1$ and $2$ together imply that $f$ sends Maurer-Cartan elements to Maurer-Cartan elements. Since $f$ respect Lie bracket, whenever $v_1,v_2\in MC(L)$ are gauge equivalent via $exp(\alpha)$ for some $\alpha\in L^0$, $f(v_1)$ and $f(v_2)$ are gauge equivalent via $exp(f(\alpha))$.
So the induced morphism on Maurer-Cartan elements descends to the quotient by gauge equivalence.
\begin{lemma}
An $L_{\infty}$-morphism $f:L\rightarrow L'$ between DGLAs induces a well defined map $Def(f):Def(L)\rightarrow Der(L')$. When $f$ is $L_{\infty}$ isomorphic ($i.e.$ $f$ has an $L_{\infty}$ inverse $g:L'\rightarrow L$) then $Def(L)\cong Def(L')$.
\begin{proof}
Notice that $f$ induces an $L_{\infty}$-morphism between $\hbar L[[\hbar]]$ and $\hbar L'[[\hbar]]$, apply above discussion.\\
When $f$ is isomorphic, the induced morphism between $\hbar L[[\hbar]]$ and $\hbar L'[[\hbar]]$ is also isomorphic. Then they have isomorphic Maurer-Cartan elements and isomorphic gauge groups (and their actions), which gives an isomorphism between $Def(L)$ and $Def(L')$.
\end{proof}
\end{lemma}
\begin{remark}
One can similarly define \textit{generalized Maurer-Cartan} elements on an $L_{\infty}$-algebra $L$. To be precise, it is the set
\[\{s\in L^1\text{ }|\text{ }l_1(s)+\frac{1}{2!}l_2(s,s)+\frac{1}{3!}l_3(s,s,s)+\cdots\}.\] A corresponding gauge equivalence can also be defined.
\end{remark}

\begin{definition}
A morphism $f=(f_1,f_2,\cdots):L\rightarrow L'$ between $L_{\infty}$-algebras is a \textit{weak equivalence} if the chain-map $f_1$ is a quasi-isomorphism. \\
A DGLA is \textit{formal} if it is weak equivalent to its cohomology group (regarded as a DGLA with zero differential).
\end{definition}
\begin{remark}
$T_{poly}$ is formal, as it has zero differential.
\end{remark}
\begin{definition}
An $L_{\infty}$-algebra $(L,(l_n)_{n\in\mathbb(N)^*})$ is \textit{minimal} if $l_1=0$. It is \textit{contractible} if $l_k=0$ for $k\geq 2$ and if the cohomology groups $H^{\bullet}(L,l_1)=0$.
\end{definition}
It is a direct verification that when $L$ is a contractible (resp. minimal) DGLA, then the induced DGLA $\hbar L[[\hbar]]$ is contractible (resp. minimal).\newline
We state without proving the following two results, whose details can be found in \cite{Kajiura_2006} and \cite{Kontsevich_2003}, respectively.
\begin{lemma}
A weak equivalence of minimal $L_{\infty}$-algebra is an $L_{\infty}$ isomorphism.
\end{lemma}
\begin{theorem}
(\textit{Decomposition theorem}) Every $L_{\infty}$-algebra $L$ can be decomposed as a direct sum of a minimal $L_{\infty}$ algebra and a contractible $L_{\infty}$ algebra $L=L_m\oplus L_c$.
\end{theorem}
For such a decomposition $L\cong L_m\oplus L_c$, since the inclusion $i:L_m\hookrightarrow L$ is a quasi-isomorphism (as $H^2(L_c)=0$), we have
\begin{corollary}
Each $L_{\infty}$-algebra is weakly isomorphic to a minimal one.
\end{corollary}
The inclusion has a \textit{quasi-inverse}, namely, the \textit{natural projection} $p:L\twoheadrightarrow L_m$.
\begin{lemma}
If $L$ is a contractible $DGLA$, then \[Def(L)= MC(\hbar L[[\hbar]])/exp(\hbar L^0[[\hbar]])\] is a singleton.
\begin{proof}
Since $L$ is contractible, $\tilde{L}= \hbar L[[\hbar]]$ is contractible. Since $\tilde{L}$ has zero Lie bracket,
\[MC(\tilde{L})=\{s=s_1\hbar+s_2\hbar^2+\cdots\text{ }|\text{ }ds_1=ds_2=\cdots=0\};\]
by exactness of $L$, for each $i$, $s_i=d\alpha_i$ for some $\alpha_i\in L^0$. Take \[\alpha:=-\alpha_1\hbar-\alpha_2\hbar^2-\cdots.\] The expansion $(15)$ of group action shows that $\alpha\cdot s=0$. So all elements of $MC(\tilde{L})$ are gauge equivalent to 0.
\end{proof}
\end{lemma}
\begin{lemma}
Let $L,L'$ be two DGLAs. Then $Def(L\oplus L')\cong Def(L)\times Def(L')$.
\begin{proof}
This follows from the fact that over $L\oplus L'$, the associated $(l_n\oplus l_n')_{n\in \mathbb{N}^*}$ and gauge group action are defined coordinate-wise.
\end{proof}
\end{lemma}
We now come to the central result of this section:
\begin{theorem}
Given two DGLAs $L$ and $L'$, a weak equivalence $f:L\rightarrow L'$ induces an isomorphism $Def(f):Def(L)\rightarrow Def(L')$.
\begin{proof}
Lemma 31 says that $Def(f)$ is well-defined. Assume $f$ is a weak equivalence. Decompose $L$ and $L'$ into direct sums $L_m\oplus L_c$ and $L_m'\oplus L_c'$ of minimal and contractible $L_{\infty}$-algebras, respectively. The composition
\[\hat{f}:L_m\overset{i}{\hookrightarrow}L_m\oplus L_c\overset{f}{\rightarrow}L_m'\oplus L_c\overset{p}{\twoheadrightarrow}L_m'\] where $i$ is the natural inclusion and $p$ is the natural projection, gives a weak equivalence from $L_m$ to $L_m'$, which by Lemma 34 is an isomorphism of DGLA. It then induces an isomorphism
\begin{multline*}
Def(\hat{f}):Def(L_m)\overset{Def(i)}{\longrightarrow}Def(L_m)\times Def(L_c)\\
\overset{Def(f)}{\longrightarrow}Def(L_m')\times Def(L_c)\overset{Def(p)}{\longrightarrow}Def(L_m').
\end{multline*}
Lemma 37 says that $Def(L_c)$ and $Def(L_c')$ are singletons, so $Def(i)$ and $Def(p)$ are isomorphic. Then $Def(f)$ is isomorphic.
\end{proof}
\end{theorem}
The path is quite clear now. Thanks to Theorem 39, once we have a weak equivalence between $T_{poly}$ and $D_{poly}$, it induces an isomorphism between $Def(T_{poly})$ (which governs the formal Poisson structure/gauge group) and $Def(D_{poly})$ (which governs the deformation product/gauge group). By HKR we have a quasi-isomorphism between them. The technical part is to extend the HKR map into a $L_{\infty}$ morphism.\newline
Let us state the theme of this thesis again:
\begin{theorem}
Given a smooth manifold $M$, there exists a weak equivalence $U=(U_1,U_2, \cdots)$ between its associated $T_{poly}$ and $D_{poly}$. Moreover, the first component $U_1$ of $U$ is the HKR isomorphism $\Psi$.
\end{theorem}
\newpage
\section{Main theorem on $\mathbb{R}^d$}
To simplify our notation, we assume \textit{Einstein summation convention}. We write $\partial_i$ for $\frac{\partial}{\partial x_i}$. For a set $S$, we write $\#S$ to be the number of elements in $S$.
\newline
On $\mathbb{R}^d$, the condition that $U=(U_1,U_2,\cdots):T_{poly}\rightarrow D_{poly}$ being $L_{\infty}$ can be written explicitly as (Please find \cite{Arnal2000ChoixDS} for a detailed analysis of signs):
\begin{multline}
    \sum_{i\neq j}\pm(U_{n-1}((\xi_i\bullet\xi_j)\wedge\xi_1\wedge\cdots\wedge\hat{\xi_i}\wedge\cdots\wedge\hat{\xi_j}\wedge\cdots\xi_n))(f_1\otimes\cdots\otimes f_m)\\
+\sum_{\substack{k,l\geq 0\\k+l=n}}\frac{1}{k!l!}\sum_{\sigma\in S_n}\pm(U_k(\xi_{\sigma(1)}\wedge\cdots\wedge\xi_{\sigma(k)})\circ U_l(\xi_{\sigma(k+1)}\wedge\cdots\wedge\xi_{\sigma(n)}))(f_1\otimes\cdots\otimes f_m)=0,
\tag{F}
\end{multline}
where $\hat{}$ stands for omission.
\subsection{Graphic representation of differential operators}
For a smooth manifold $M$ with $A:=C^{\infty}(M)$,
We want to assign polyvector fields on $M$ with linear polydifferential operators on $A$. As an example, for $\xi_1\in\bigwedge^2Der(A),\xi_2\in\bigwedge^3Der(A)$, we can consider the operator
\begin{align}
    (f_1,f_2,f_3)\mapsto \xi_2^{i^2_1i^2_2}\partial_{i^2_1}\xi_1^{i^1_1i^1_2i^1_3}\partial_{i^1_1}f_1\partial_{i^1_2}f_2\partial_{i^1_3}\partial_{i^2_2}f_3
\end{align}
This sort of polydifferential operators can be reformulated by graphs. A graphic representation of the above example looks like:
\begin{figure}[htp]
    \centering
    \includegraphics[width=8cm]{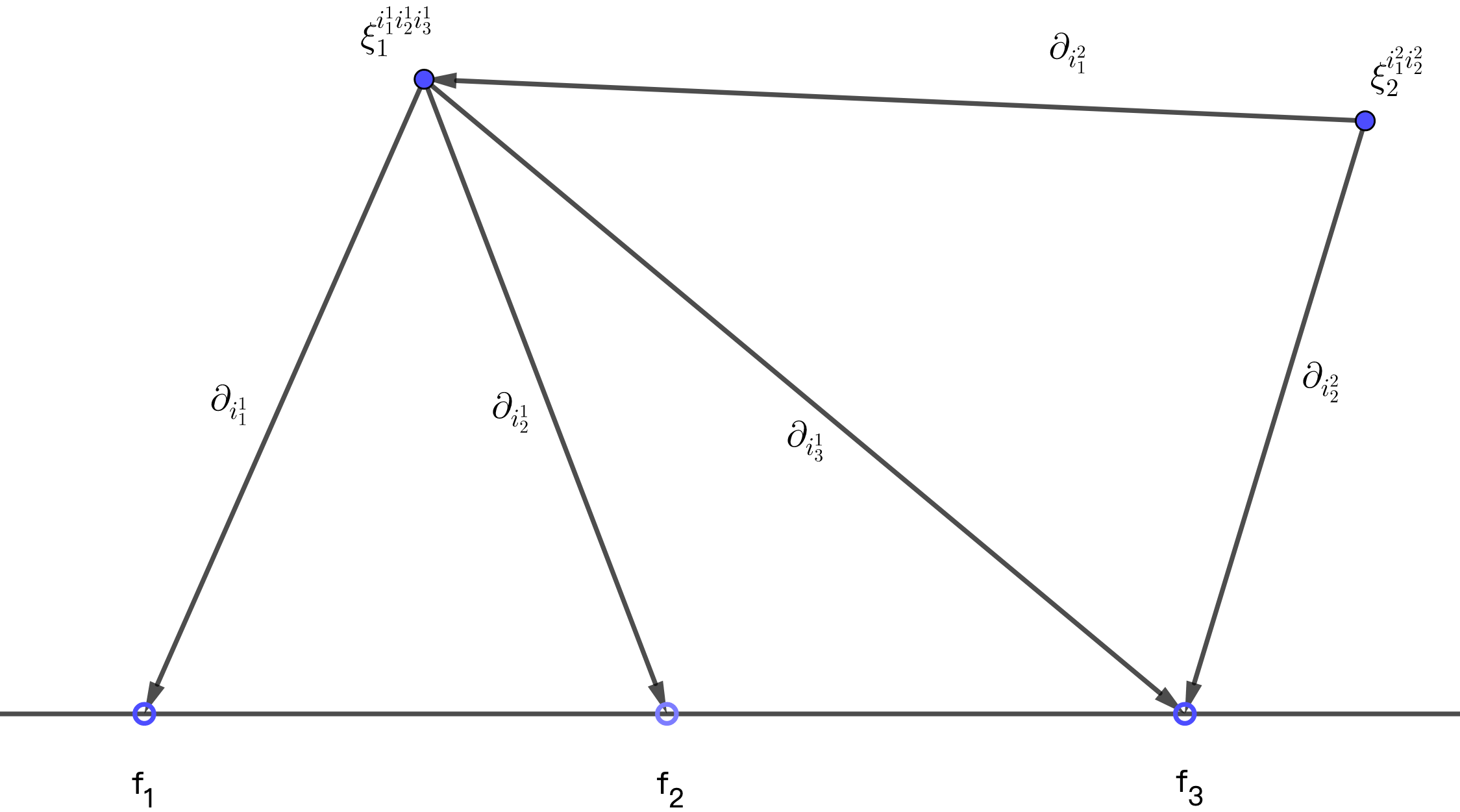}
    \caption{A graphic representation of (16)}
    \label{fig:graphic representation example}
\end{figure}
\newline
And a graph like the above assigns a pair $(\xi_1\in\bigwedge^3(Der(A)),\xi_2\in\bigwedge^2(Der(A)))$ with an element of $D^3(A,A)$.
\newline
In fact, a \textit{suitably} defined graph assigns a linear polydifferential operator to an element of $T_{poly}$. Such graphs are said to be \textit{admissible}.
\begin{figure}[htp]
    \centering
    \includegraphics[width=8cm]{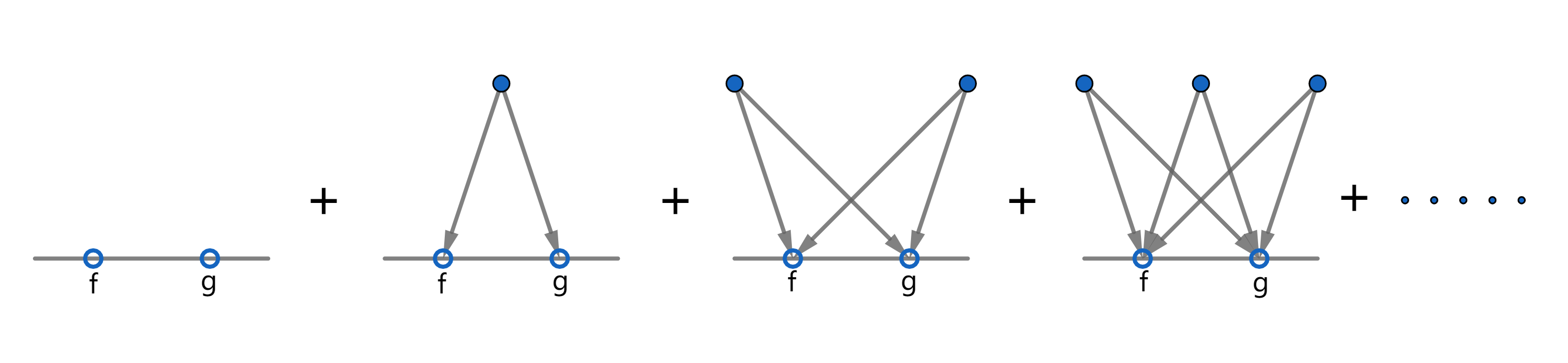}
    \caption{The Moyal $f*g$ product can be represented by graphs as shown. Each solid vertex represents the constant Poisson bivector field $\Pi$, and each outgoing arrow represents a derivative.}
    \label{fig:Moyal graphic representation}
\end{figure}
\begin{definition}
An oriented graph $\Gamma$ is \textit{admissible}, if it satisfies the following properties:
\begin{enumerate}
    \item The set of vertices $V(\Gamma)$ is divided into two ordered subsets: $V_1(\Gamma)$ of order $n$, and $V_2(\Gamma)$ of order $\bar{n}$, called vertices of \textit{first type} and \textit{second type} respectively.
    \item The orders $n$, $\bar{n}$ satisfy $2n+\bar{n}-2\geq 0$.
    \item The set of edges $E(\Gamma)$ is finite.
    \item For each edge in $E(\Gamma)$, the starting point and ending point do not coincide.
    \item All edges start from vertices of first type.
    \item The set $star(v)$ of edges coming out from a given vertex $v\in V_1(\Gamma)$ is ordered.
\end{enumerate}
Sometimes we write $\Gamma_n$ (resp. $\Gamma_{\bar{n}}$) for an admissible graph with n vertices of first type (resp. $\bar{n}$ vertices of second type). We denote $G_{n,\bar{n}}$ for the set of all admissible graphs with $n$ vertices of first type and $\bar{n}$ vertices of second type.
\end{definition}
With $\Gamma\in G_{n,\bar{n}}$, we can assign a polydifferential operator to a tensor product of (homogeneous) polyvector fields \[\xi_1\otimes\cdots\otimes\xi_n\in\bigwedge\nolimits^{\#(star(v_1))}Der(A)\otimes\cdots\otimes\bigwedge\nolimits^{\#(star(v_n))}Der(A)\subset (\bigwedge\nolimits^{\bullet}Der(A))^{\otimes n}\] via the following procedure:
\begin{enumerate}
    \item Associate to each vertex $v_j$ of first type ($j=1,\cdots,n$) with $k=\#(star(v_j))$ outgoing arrows the skew-symmetric tensor $\xi_j^{i^j_1i^j_2\cdots i^j_k}$.
    \item For the $l^{th}$ edge of $star(v_j)$, associate a partial derivative $\partial_{i^j_l}$.
    \item Place a function $f_{\bar{j}}$ to each vertex $v_{\bar{j}}$ of second type ($\bar{j}=1,\cdots,\bar{n}$).
    \item Multiply such elements in the order prescribed by the labelling of the graph.
\end{enumerate}
The procedure associates, for each $\Gamma\in G_{n,\bar{n}}$, a linear map
\[U_{\Gamma}':\bigwedge\nolimits^{\#(star(v_1))}Der(A)\otimes\cdots\otimes\bigwedge\nolimits^{\#(star(v_n))}Der(A)\rightarrow D^{\bar{n}}(A),\]
which can be extended (by zero) to a homogeneous, multilinear map
\[U_{\Gamma}:(\bigwedge\nolimits^{\bullet}Der(A))^{\otimes n}\rightarrow D^{\bar{n}}(A).\]
A permutation $\sigma\in S_n$ acts on $\Gamma$ by permuting the order of vertices of first type. Anti-symmetrizing $U_{\Gamma}$ w.r.t the group action of $S_n$ gives rise to a (graded) skew-symmetric multilinear map $\tilde{U}_{\Gamma}$, so that the expression $\tilde{U}_{\Gamma}(\xi_1\wedge\cdots\wedge\xi_n)$ makes sense. But we are not going to do so at this moment.
\newline
The desired $L_{\infty}$ morphism $U=(U_1,U_2,\cdots)$ has its $n^{th}$ component of the form
          \[U_n:=\sum_{\bar{n}\geq 0}\sum_{\substack{\Gamma\in G_{n,\bar{n}}\\\#(E(\Gamma))=2n+\bar{n}-2}}W_{\Gamma}\times U_{\Gamma},\]
for some $W_{\Gamma}$ to be determined. The condition for $U_n$ being graded skew-symmetric is encoded in the construction of $W_{\Gamma}$\footnote{That's why we don't anti-symmetrize at first.}, as we will see. We will also show that $U_n$ has the desired degree.
\subsection{Configuration spaces}
\subsubsection{Definitions}
Let $n,\bar{n}$ be non-negative integers satisfying $2n+\bar{n}\geq 2$. By $Conf_{n,\bar{n}}$ we mean the product of the \textit{configuration space of the upper half plane} with the \textit{configuration space of the real line}:
\begin{align*}
Conf_{n,\bar{n}}:=\{(p_1,\cdots,p_n;q_1,\cdots,q_{\bar{n}})|p_j\in\mathbb{H}_{+}^2, q_j\in\mathbb{R},p_{j_1}\neq p_{j_2} \text{ for }j_1\neq j_2, q_{\bar{j}_1}\neq q_{\bar{j}_2} \text{ for } \bar{j}_1\neq \bar{j}_2\}.
\end{align*}
Here $\mathbb{H}^2_{+}\subset\mathbb{R}^2$ stands for the upper half plane in $\mathbb{R}^2$ without $x$-axis, identified with $\{x+iy\in\mathbb{C}| y>0\}$.\newline
$Conf_{n,\bar{n}}$ is a smooth manifold of real-dimension $2n+\bar{n}$. We consider points in $Conf_{n,\bar{n}}$ up to \textit{rescaling} and \textit{translation}. More precisely, we consider the Lie group $G:=\mathbb{R}^+\times \mathbb{R}$ acting on $Conf_{n,\bar{n}}$ defined as
\[(a,b)\cdot z:=az+b\]
for $(a,b)\in G$ and $z\in Conf_{n,\bar{n}}$. \newline
The action of $G$ on $Conf_{n,\bar{n}}$ is free, since there is no point in $\mathbb{H}^2_{+}$ that is invariant under rescaling or translation. \\
As $G$ is a Lie group of real dimension $2$ acting freely on $Conf_{n,\bar{n}}$ of real dimension $2n+\bar{n}>2$, the quotient \[C_{n,\bar{n}}:=Conf_{n,\bar{n}}/G\]
is again a smooth manifold, whose dimension is $2n+\bar{n}-2$. \newline
We can similarly define
\[Conf_n:=\{(p_1,\cdots,p_n)|p_i\in\mathbb{C}, p_{j_1}\neq p_{j_2}\text{ for }j_1\neq j_2\},\]
a smooth manifold of real dimension $2n>4$ (since $n>2$), on which the 3-dimensional-Lie group $G'=\{\mathbb{R}^+\times\mathbb{C}\}$ acts freely, via
\[(a,b)\cdot z:=az+b.\]
We can define analogously the quotient space $C_n:=Conf_n/G'$, but this time $C_n$ has real dimension $2n-3$.\\
Finally, we let $C^+_{n,\bar{n}}$ to be the submanifold of $C_{n,\bar{n}}$ consisting of points \[(p_1,\cdots,p_n;q_1,\cdots,q_{\bar{n}})\] satisfying $q_1<\cdots <q_{\bar{n}}$. Since points of this sort in $Conf_{n,\bar{n}}$ form an open connected submanifold 
of $Conf_{n,\bar{n}}$ and the action of $G$ does not affect the order of $\{q_{\bar{j}}\}_{1\leq \bar{j}\leq \bar{n}}$, we conclude that $C_{n,\bar{n}}^+$ is connected.\newline
There is a natural correspondence of vertices in an admissible graph $\Gamma\in G_{n,\bar{n}}$ to coordinates points of $C^+_{n,\bar{n}}$: simply by identifying the ordered sets $V_1(\Gamma)$ and $V_2(\Gamma)$ with $\{p_1,\cdots,p_n\}$ and $\{q_1,\cdots,q_{\bar{n}}\}$, respectively.

\subsubsection{Compactification and boundary strata}
We are going to define $W_{\Gamma}$ via integration of differential forms on compactifications of $C_{n,\bar{n}}$ and $C_n$. We state the following results about compactifying these spaces. Please find \cite{Axelrod:1993wr} for more details.
\begin{theorem}
For any $n,\bar{n}\in\mathbb{N}$ satisfying $n+\bar{n}>2$, there exists a compactification $\bar{C}_{n,\bar{n}}$ (resp. $\bar{C}_n$ where $n>2$) whose interior is the open configuration space $C_{n,\bar{n}}$ (resp. $C_n$).
\end{theorem}
\begin{remark}
The compactification of configuration spaces is at first introduced in the algebro-geometric setting in the paper \cite{Fulton:1994hh}.
\end{remark}
The process of compactification of $C_{n,\bar{n}}$ and $C_n$ (as smooth manifolds) is basically done by embedding them into some larger manifolds followed by taking compactification of the image:
\newline
\textit{The case $C_n$ where $n>2$}: given any point $(p_1,\cdots,p_n)$ of $C_n$ we associate a collection of $n(n-1)$ angles
\[\{Arg(p_i-p_j)|i\neq j\}\]
whose values are in $\mathbb{R}/2\pi\mathbb{Z}\cong S^1$;\\
and a collection of $n^2(n-1)^2$ ratios of distances
\[\{\frac{||p_i-p_j||}{||p_k-p_l||} | i\neq j,k\neq l\}\]
whose values are in $[0,+\infty)$. \\
With these angles and ratios we obtain an embedding of $C_n$ into the manifold $(S^1)^{n(n-1)}\times [0,+\infty)^{n^2(n-1)^2}$.\newline
\textit{The case $C_{n,\bar{n}}$}: we first embed $C_{c,\bar{n}}$ into $C_{2n+\bar{n}}$ by the mapping on the level of configuration spaces:
\[Conf_{n,\bar{n}}\ni (p_1,\cdots,p_n;q_1,\cdots, q_{\bar{n}})\mapsto (p_1,\cdots,p_n,\bar{p_1},\cdots,\bar{p_n};q_1,;q_1,\cdots, q_{\bar{n}})\in Conf_{2n,\bar{n}},\]
and $\bar{C}_{n,\bar{n}}$ is (defined to be) the compactification of the image in $C_{2n,\bar{n}}$.
\begin{definition}
A \textit{manifold with corners of dimension m} is a Hausdorff topological space $M$ which is locally homeomorphic to $\mathbb{R}^{m-n}\times [0,\infty)^n$ where $n=0,1,\cdots, m$. If the local homeomorphisms are diffeomorphic, $M$ is said to be smooth.\\
A point $x$ in a manifold with corners $M$ of dimension $m$ is said to be of type $n$ if $x=(x_1,x_2,\cdots,x_{m-n},0,\cdots,0)$ in some (hence any) chart. All points in $M$ of type $n$ form a submanifold of $M$ called \textit{stratum of codimension $n$}. In particular, the stratum of codimension $1$ in $M$ is the boundary $\partial M$ of $M$.
\end{definition}
The compactifications $C_{n,\bar{n}}$ and $C_n$ are smooth manifolds with corners. Intuitively, we can think of $\partial \bar{C}_{n,\bar{n}}$ as \textit{degeneration} of $C_{n,\bar{n}}$ by means of \textit{collapsing} some of the coordinates points. This gives a classification of $\partial \bar{C}^+_{n,\bar{n}}$ into two types:
\begin{enumerate}
    \item \textit{Strata of type S1}: Points $p_j\in\mathbb{H}^2$ for $j\in S\subset\{1,\cdots,n\}$, where $\#S>2$, collapse into a point in $\mathbb{H}^2$. Strata of this type (locally) have the form
    \[C_{\#S}\times C_{n-\#S+1,\bar{n}},\]
The first term represents the relative position of the collapsing points from $S$, when we are looking at them through a magnifying glass. The second term is the space of all remaining points plus the new point formed by collapsing.
\begin{figure}[h]
    \centering
    \includegraphics[width=8cm]{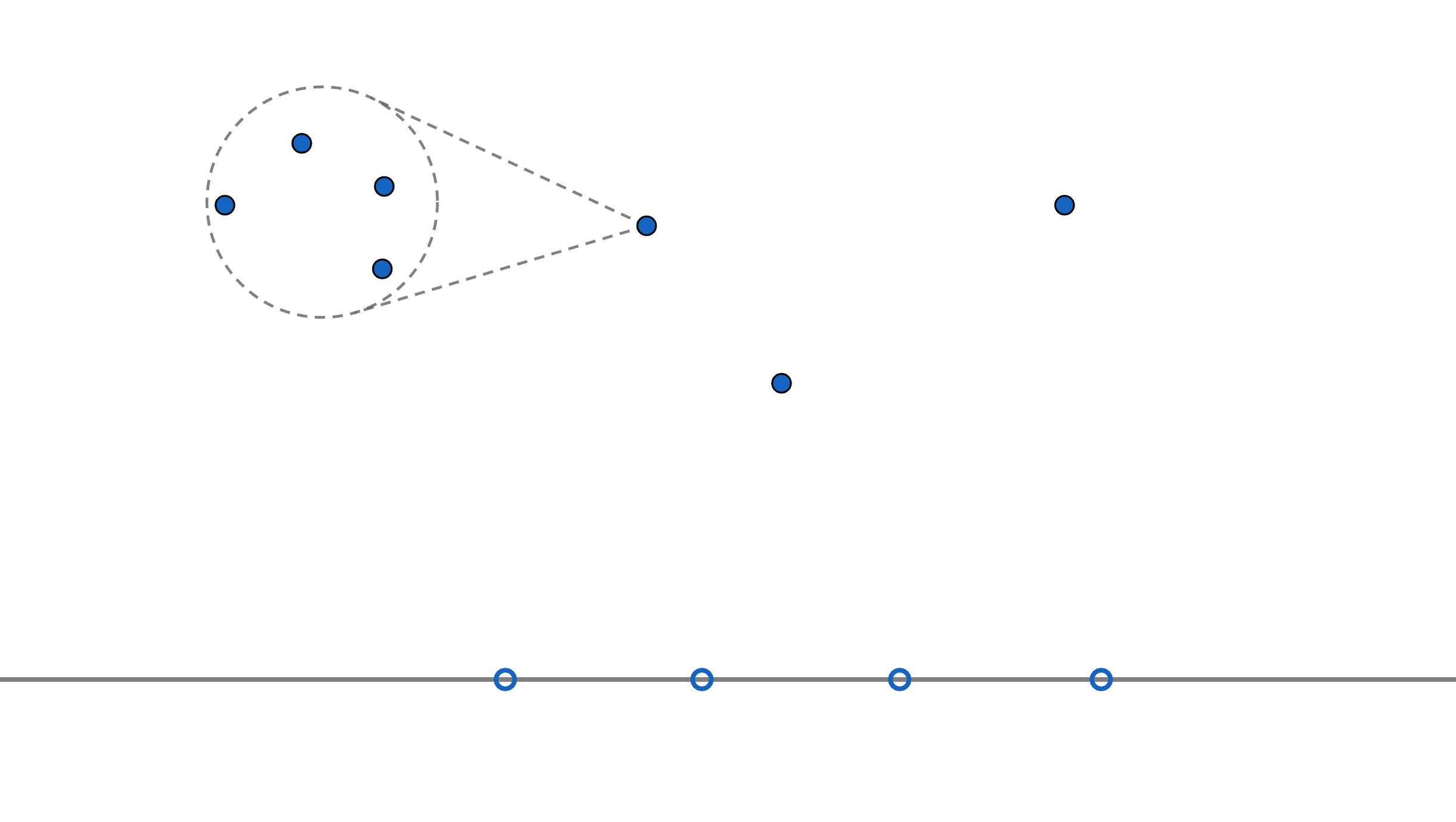}
    \caption{A type S1 stratum}
    \label{ }
\end{figure}
\newline
On an admissible graph $\Gamma\in G_{n,\bar{n}}$, this is to contract a subgraph $\Gamma'\in G_{\#S,0}$ with vertices from $S\subset V_1(\Gamma)$ (via the identification) into a new vertex of first type. Under this interpretation, the first term corresponds\footnote{Notice here identifying with coordinates of $C_{\#S}$ instead of $C_{\#S,0}$} to the data about the contracted subgraph $\Gamma'$, while the second term corresponds to the data about the newly-formed graph $\tilde{\Gamma}\in G_{n-\#S+1,\bar{n}}$ after contraction.
\begin{figure}[h]
    \centering
    \includegraphics[width=8cm]{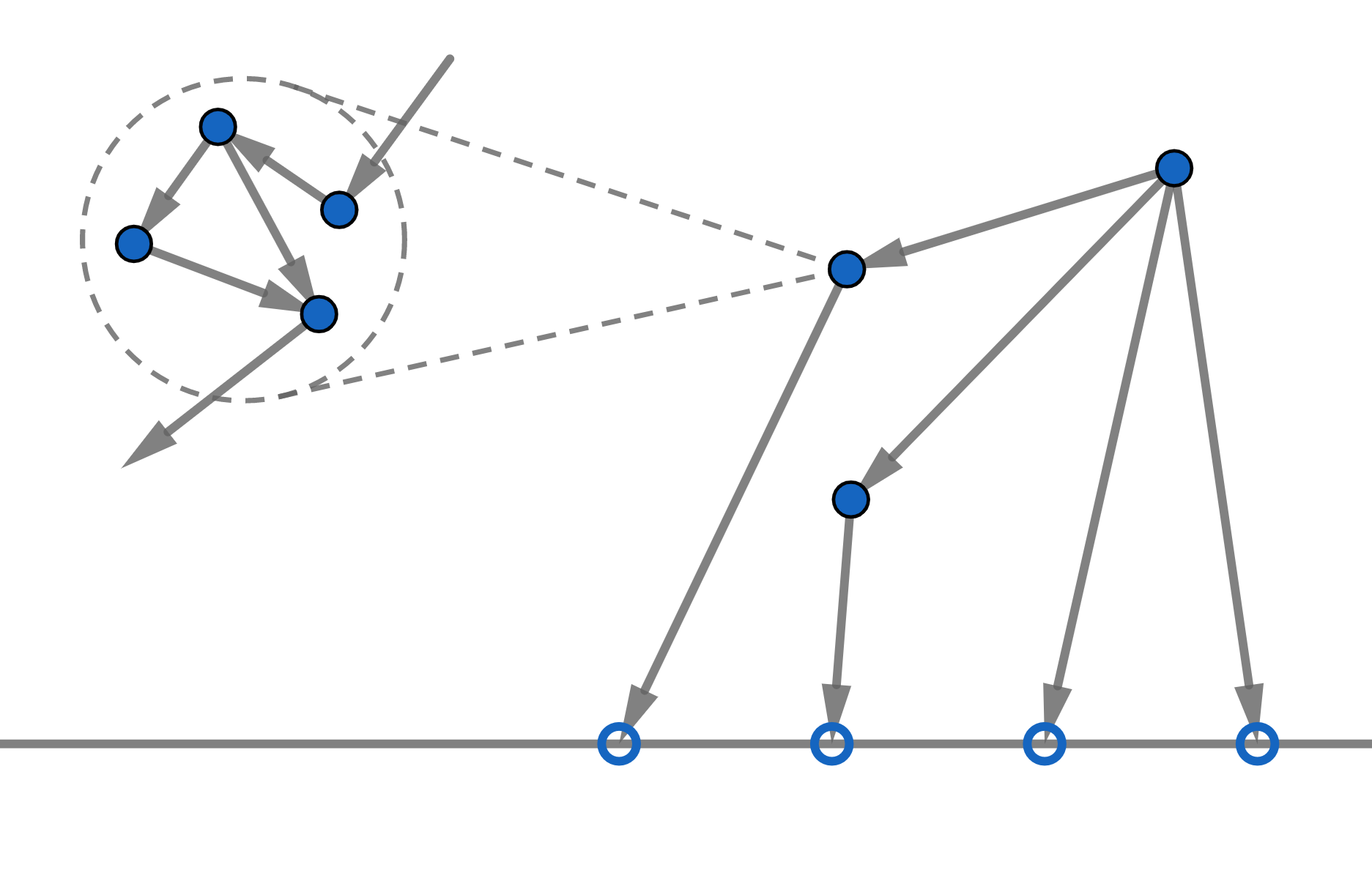}
    \caption{Contraction of a graph corresponding to a type S1 stratum}
    \label{ }
\end{figure}
\newline
\item \textit{Strata of type S2}: Points $p_j\in\mathbb{H}^2$ for $j\in S\subset\{1,\cdots,n\}$ and points $q_{\bar{j}}\in\mathbb{R}$ for $\bar{j}\in\bar{S}\subset\{1,\cdots,\bar{n}\}$, where $2\#S+\#\bar{S}\geq 2$,
collapse into a point in $\mathbb{R}$. strata of this type (locally) have the form
\[C_{\#S,\#\bar{S}}\times C_{n-\#S,\bar{n}-\#\bar{S}+1}\]
\begin{figure}[htp]
    \centering
    \includegraphics[width=8cm]{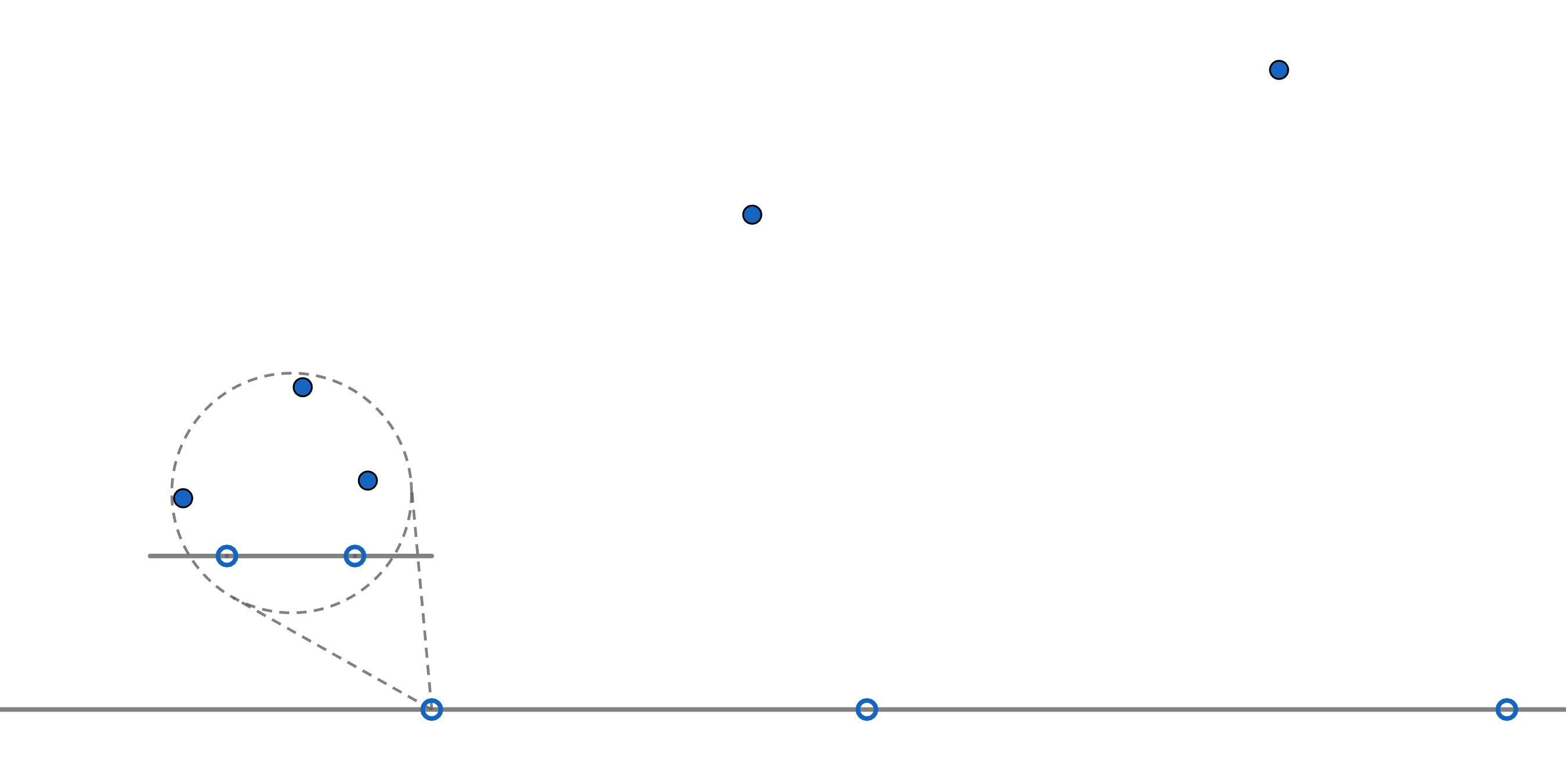}
    \caption{A type S2 stratum}
    \label{ }
\end{figure}
On an admissible graph $\Gamma\in G_{n,\bar{n}}$, this is to contract a subgraph $\Gamma'\in G_{\#S,\#\bar{S}}$ with vertices of first type from $S\subset V_1(\Gamma)$ and vertices of second type from $\bar{S}\in V_2(\Gamma)$ into a new vertex of second type. The first term corresponds to the data about the contracted subgraph $\Gamma'$, while the second corresponds to the data about the newly-formed graph $\tilde{\Gamma}\in G_{n-\#S,\bar{n}-\#\bar{S}+1}$ after contraction.
\end{enumerate}
%
%
%
\begin{figure}[htp]
    \centering
    \includegraphics[width=8cm]{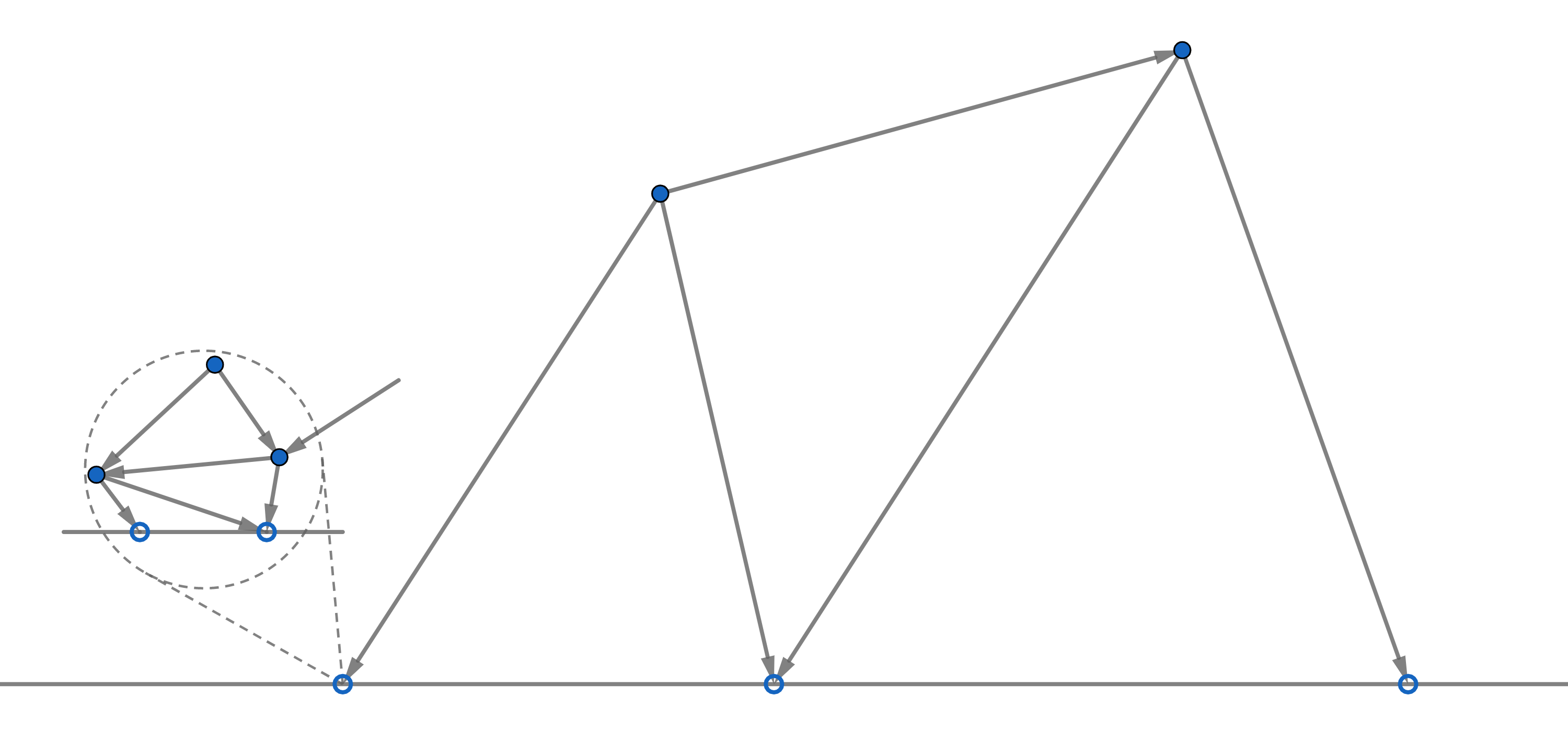}
    \caption{Contraction of a graph corresponding to a type S2 stratum}
    \label{ }
\end{figure}
We can look at examples of $C_{n,\bar{n}}$ and $C_n$ and their compactifications for small $n$'s and $\bar{n}$'s.\\
\textbf{Examples}:
\begin{enumerate}
    \item $C_{1,0}\cong \mathbb{H}^2/G$ is a single point, since, for any $p\in\mathbb{H}^2$, there exists an affine transformation moving $p$ to the point $i=\sqrt{-1}$. $\bar{C}_{1,0}=C_{1,0}$ is also a single point.
    \item $C_{0,2}$ consists of two points (on the real line $\mathbb{R}$). Hence $\bar{C}_{0,2}=C_{0,2}$ is also a two-point-space.
    \item $C_{1,1}$ is (homeomorphic to) an open interval: for any $(p,q)\in Conf_{1,1}$ there is a unique affine transformation moving $(p,q)$ to a point on $\{(e^{i\pi\theta},0)|0<\theta<1\}$. $\bar{C}_{1,1}$ is then (homeomorphic to) a closed interval.
    \item $C_2$ is homeomorphic to $S^1$: For any $(p_1,p_2)\in Conf_{2}$ there is a unique element in $G'$ moving $(p_1,p_2)$ to $(\frac{p_1-p_2}{||p_1-p_2||},0)$.
    \item $C_{2,0}\cong \mathbb{H}^2-\{(0,i)\}$: We have for every $(p_1,p_2)\in Conf_{2,0}$ a unique affine transformation moving $p_2$ to $i$, then $p_1$ is still a point on $\mathbb{H}^2$ but $p_1\neq p_2$.\\
    By contracting coordinate points of $C_{2,0}$ we can explicitly decompose its compactification:
    \[\bar{C}_{2,0}=C_{2,0}\sqcup C_{1,1}\sqcup C_{1,1}\sqcup C_{0,2}\sqcup C_2.\]
    In particular, $\bar{C}_{1,1}=C_{0,2}\sqcup C_{1,1}\subset \bar{C}_{2,0}$.
    %
    %
    \begin{figure}[htp]
    \centering
    \includegraphics[width=8cm]{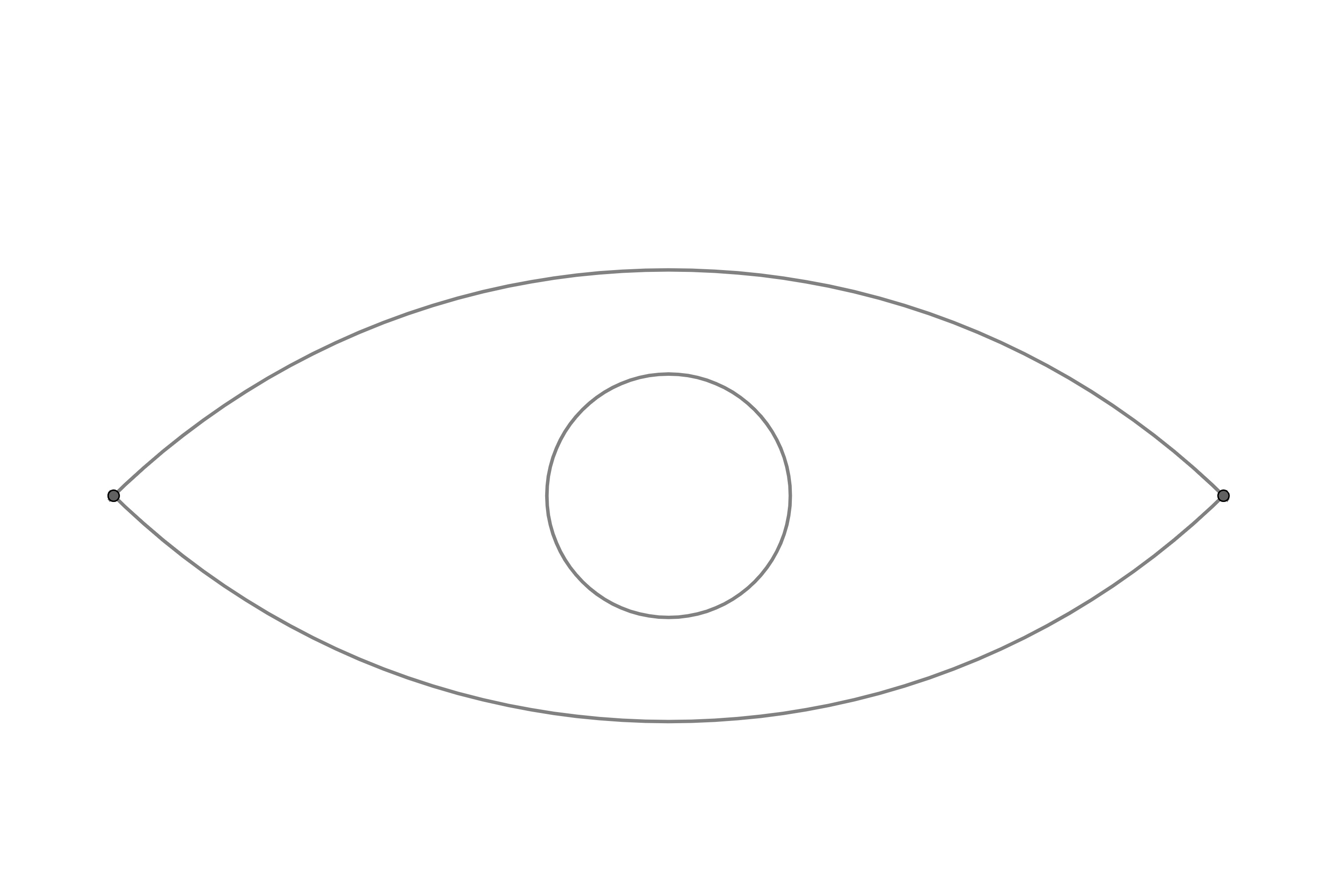}
    \caption{Compactification of $C_{2,0}$}
    \label{fig:Kontsevich's eye}
\end{figure}
\end{enumerate}
\subsubsection{Differential forms on $\bar{C}_{n,\bar{n}}$}
To each pair of distinct points $p_1,p_2\in\mathbb{H}^2$ we can define an \textit{angle map}
\[\phi:(p_1,p_2)\mapsto \frac{1}{2i}Log(\frac{(p_2-p_1)(\bar{p_2}-p_1)}{(p_2-\bar{p_1})(\bar{p_2}-\bar{p_1})})=Arg(\frac{p_2-p_1}{p_2-\bar{p_1}})\in \mathbb{R}/2\pi\mathbb{Z}\cong S^1.\]
Geometrically, this is to equip $\mathbb{H}^2$ with the \textit{Poincaré metric}, on which the geodesic line between $p_1,p_2$ is a circular arc connecting $p_1,p_2$ where the circle has its center at $x$-axis. now draw a vertical line through $p_1$ and a geodesic line from $p_1$ to $p_2$, then $\phi(p_1,p_2)$ is the angle between the vertical line and the tangent line of the geodesic line at $p_1$ in the usual sense. See below:
%
%
\begin{figure}[htp]
    \centering
    \includegraphics[width=8cm]{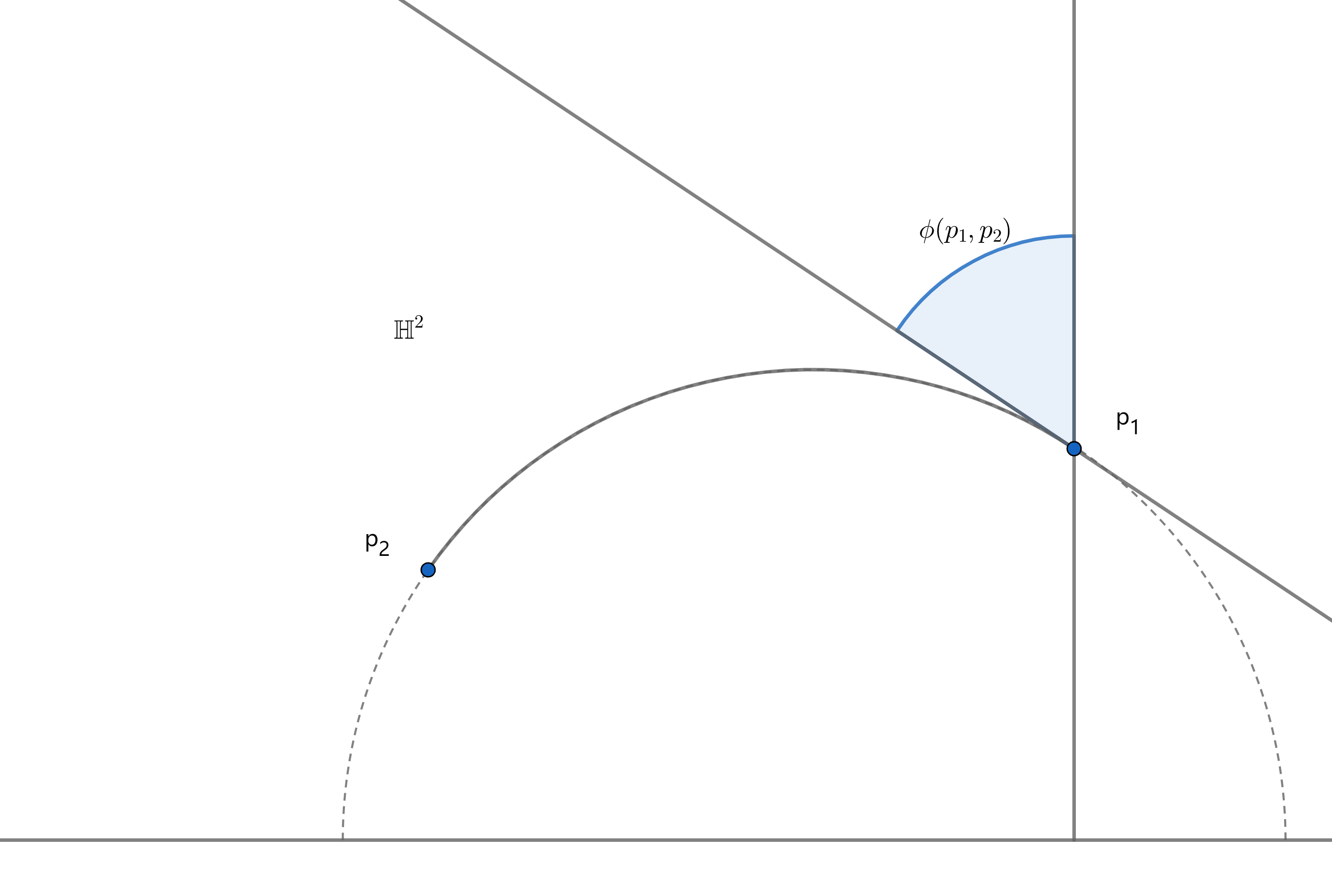}
    \caption{An illustration for $\phi(p_1,p_2)$ }
    \label{fig: angle measure}
\end{figure}
\begin{remark}
When $p_1$ lies on $x$-axis, from the illustration we have $\phi(p_1,p_2)$ is zero.
\end{remark}
We thus obtain a smooth map from $C_{2,0}$ to $S^1$, which then extends smoothly \footnote{See \cite{Kontsevich_2003}}
to a map (still call it $\phi$) from $\bar{C}_{2,0}$ to $S^1$.
The differential $d\phi$ is a smooth $1$-form on $\bar{C}_{2,0}$.
\newline
On the other hand, by the identification between vertices of an admissible graph $\Gamma\in G_{n,\bar{n}}$ and coordinate points of $C_{n,\bar{n}}$, every oriented edge $e$ assigns a projection from $C_{n,\bar{n}}$ to $\bar{C}_{2,0}$ (which then extend smoothly to their boundaries):
\begin{enumerate}
    \item If $e$ is an arrow $v_{j_1}\rightarrow v_{j_2}$ between points of first type, we associate the projection $p_e:C_{n,\bar{n}}\rightarrow C_{2,0}$ by
    \[(p_1,\cdots,p_n;q_1,\cdots,q_{\bar{n}})\mapsto (p_{j_1},p_{j_2})\]
    then compose with $C_{2,0}\hookrightarrow \bar{C}_{2,0}$ (still call it $p_e$).
    \item If $e$ is an arrow $v_j\rightarrow u_{\bar{j}}$ from a point of first type to a point of second type, we associate the projection $p_e: C_{n,\bar{n}}\rightarrow C_{2,0}$ by
    \[(p_1,\cdots,p_n;q_1,\cdots,q_{\bar{n}})\mapsto (p_j,q_{\bar{j}})\in C_{1,1},\]
    then compose with $C_{1,1}\hookrightarrow \bar{C}_{2,0}$ (still call it $p_e$).
\end{enumerate}
For each edge, the differential of the composition $\phi_e:\bar{C}_{n,\bar{n}}\overset{p_e}{\longrightarrow} \bar{C}_{2,0}\overset{\phi}{\longrightarrow} S^1$ is a smooth $1$-form on $\bar{C}_{n,\bar{n}}$. We let $w_{\Gamma}$ be the differential form
\[w_{\Gamma}:=\bigwedge_{e\in E(\Gamma)}d\phi_e.\]
\begin{remark}
By Leibniz rule, $w_{\Gamma}$ is a closed form. Also, when $\#E(\Gamma)=2n+\bar{n}-2$, $w_{\Gamma}$ is a top-form on $C_{n,\bar{n}}$.
\end{remark}
The coefficient $W_{\Gamma}$ for $U_{\Gamma}$ in $U_n$ is defined as
\begin{align}
W_{\Gamma}:=\prod_{j=1}^n\frac{1}{(\#star(v_j))!}\frac{1}{(2\pi)^{2n+\bar{n}-2}}\int_{\bar{C}^+_{n,\bar{n}}}w_{\Gamma}.
\end{align}
$U_n$ is graded skew-symmetric: if a graph $\Gamma'$ is obtained from $\Gamma$ by swapping two points $v_{j_1},v_{j_2}$ in $V_1(\Gamma)$, then $w_{\Gamma'}=-(-1)^{(\#star(v_{j_1}))(\#star(v_{j_2}))}w_{\Gamma}$ by nature of wedge product.\\
So it makes sense to write $U_n(\xi_1\wedge\cdots\wedge\xi_n)$ for polyvector fields
$\xi_1\cdots,\xi_n\in\bigwedge\nolimits^{\bullet}Der(A)$.
\subsection{Proof of the main theorem on $\mathbb{R}^d$}
To show the sequence $(U_n)_{n\geq 1}$ is an $L_{\infty}$ quasi-isomorphism, we need to verify the following:
\begin{enumerate}
    \item each $U_n$ has degree $1-n$;
    \item $U_1$ is a quasi-isomorphism (in fact, it is the HKR quasi-isomorphism $\Psi$ appears in section $4$.)
    \item $(U)_{n\geq 1}$ satisfies equations $(F)$ appear in the beginning of this section.
\end{enumerate}
\begin{lemma}
Each $U_n$ has degree $1-n$.
\begin{proof}
Recall that $U_n$ is a linear combination of $U_{\Gamma}$, where $\Gamma\in G_{n,\bar{n}}$ for some $\bar{n}$ satisfying $\sum\limits^{n}_{j=1}\#star(v_j)=E(\Gamma)=2n+\bar{n}-2$, which implies $\bar{n}-\#E(\Gamma)=2-2n$. On the other hand, the nonzero part of $U_{\Gamma}$ takes in a tensor of elements $\xi_1\otimes\cdots\otimes\xi_n$ whose degree is $\sum\limits^{n}_{j=1}(\#star(v_j)-1)=\#E(\Gamma)-n$ in $T_{poly}$, giving an element in $D^{\bar{n}}(A)$, whose degree is $\bar{n}-1$ in $D_{poly}$. So the degree of $U_{\Gamma}$ is $\bar{n}-1-\#E(\Gamma)+n=2-2n-1+n=1-n$.
\end{proof}
\end{lemma}
\begin{lemma}
$U_1=\Psi$, where $\Psi$ is the HKR quasi-isomorphism given by
\begin{align*}
\Psi_n:\xi_0\wedge\cdots\wedge\xi_n\mapsto[(f_0\otimes\cdots\otimes f_n)\mapsto\frac{1}{(n+1)!}\sum_{\sigma\in S_{n+1}}sgn(\sigma)\prod_{i=0}^{n}\xi_{\sigma(i)}(f_i)],
\end{align*}
\begin{proof}
$U_1$ is given by admissible graphs $\Gamma_{\bar{n}}$ with $n=1$ and $\#E(\Gamma_{\bar{n}})=\bar{n}$ (see below).\newline
%
%
\begin{figure}[h]
    \centering
    \includegraphics[width=8cm]{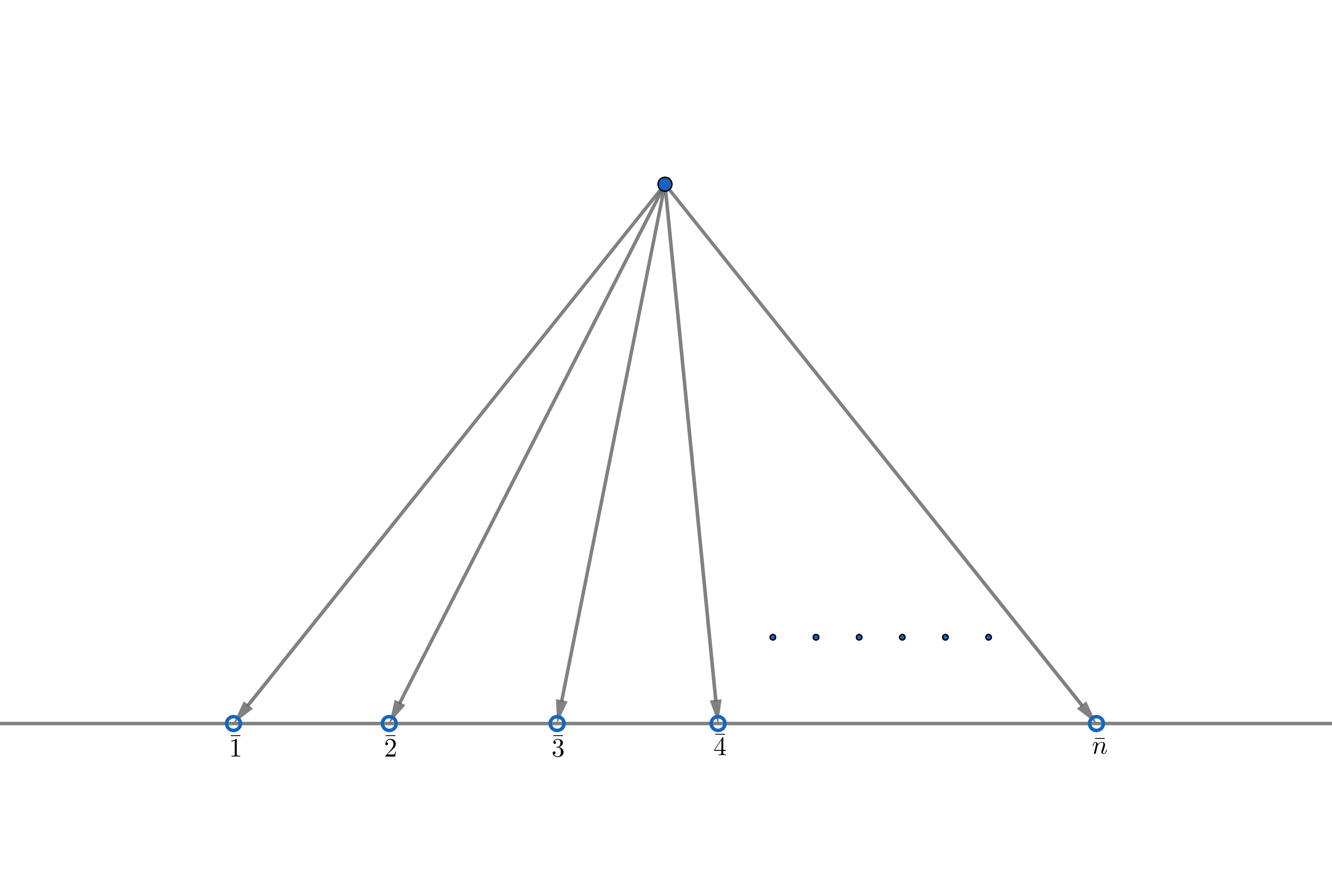}
    \caption{The only graph $\Gamma_{\bar{n}}\in G_{1,\bar{n}}$ }
    \label{fig: U1}
\end{figure}
\newline
We have, for $\xi=\xi_1\wedge\cdots\wedge\xi_{\bar{n}}\in \bigwedge^{\bar{n}}Der(A)$,
\[U_{\Gamma_{\bar{n}}}(\xi)=\xi^{i_1\cdots i_{\bar{n}}}\partial_{i_1}f_1\cdots\partial_{i_{\bar{n}}}f_{\bar{n}}\]
which is equal to $\sum_{\sigma\in S_{n}}sgn(\sigma)\prod_{i=1}^{n}\xi_{\sigma(i)}(f_i)$ by converting Einstein summation into the usual sum.\\
We can calculate $W_{\Gamma_{\bar{n}}}$ directly:
\[W_{\Gamma_{\bar{n}}}=\frac{1}{\bar{n}!}\frac{1}{(2\pi)^{\bar{n}}}
\int_{\bar{C}_{1,\bar{n}}}\bigwedge_{e\in E(\Gamma_{\bar{n}})}d\phi_e=\frac{1}{\bar{n}!}\frac{1}{(2\pi)^{\bar{n}}}\times (2\pi)^{\bar{n}}=\frac{1}{\bar{n}!}\]
Since the integral over $\bar{C}_{1,\bar{n}}$ can be decomposed as a product of $\bar{n}$ copies of integrating $d_{\phi_e}$ over $S^1$.\\
Therefore $W_{\Gamma_{\bar{n}}}\times U_{\Gamma_{\bar{n}}}$ is indeed the $(n-1)^{th}$ component $\Psi_{n-1}$ of $\Psi$.
\end{proof}
\end{lemma}
\subsubsection{Interpretation of $(F)$}
Let us have a glance through $(F)$ before our final step of verification.\newline
Summands in the first term of $(F)$ come from contraction of single edges between points of first type in some admissible graph $\Gamma'$: suppose $v_{j_1},v_{j_2}\in V_1(\Gamma')$ together with an arrow $v_{j_1}\rightarrow v_{j_2}$ lies in the $l^{th}$ place in $star(v_{j_1})$. In corresponding $U_{\Gamma}$, the arrow contributes to a term $\frac{\partial \xi_{j_1}}{\partial \zeta_{i_l}}\frac{\partial \xi_{j_2}}{\partial x_{i_l}}$\footnote{This is NOT an Einstein notation.}, which is a summand in $\xi_{j_1}\bullet\xi_{j_2}$. Therefore, by taking the sum over all permutations of $\#star(v_{j_1})$ we can replace (both vertices and the arrow) $v_{j_1}\rightarrow v_{j_2}$ in $\Gamma'$ by a new point $v$, whose corresponding polyvector field is $\xi_{j_1}\bullet\xi_{j_2}$.
On the other hand, a graphic representation of $\xi_{j_1}\bullet\xi_{j_2}$ is to draw an arrow $v_{j_1}\rightarrow v_{j_2}$, so that $\xi_{j_1}\bullet\xi_{j_2}$ appears in the sum of resulting $U_{\Gamma'}$s over all possible permutations of $star(v_{j_1})$.\\
Such $\Gamma'$ has $n=n-1+1$ vertices of type $1$, $m$ vertices of type $2$, and $\#E(\Gamma')=2(n-1)+m-2+1=2n+m-3$.  \newline
%
%
\begin{figure}[htp]
    \centering
    \subfloat{{\includegraphics[width=5cm]{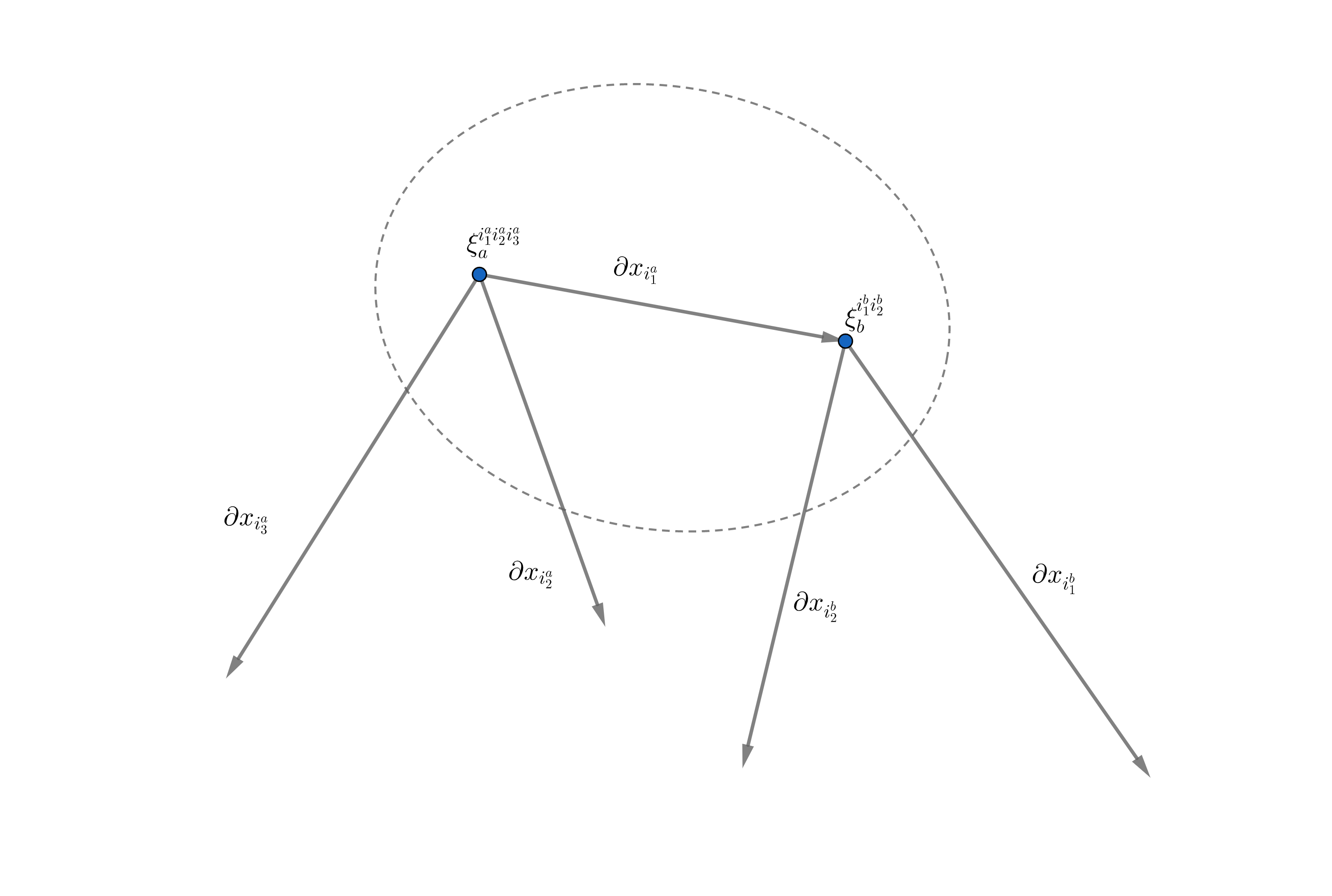} }}
    \qquad
    \subfloat{{\includegraphics[width=5cm]{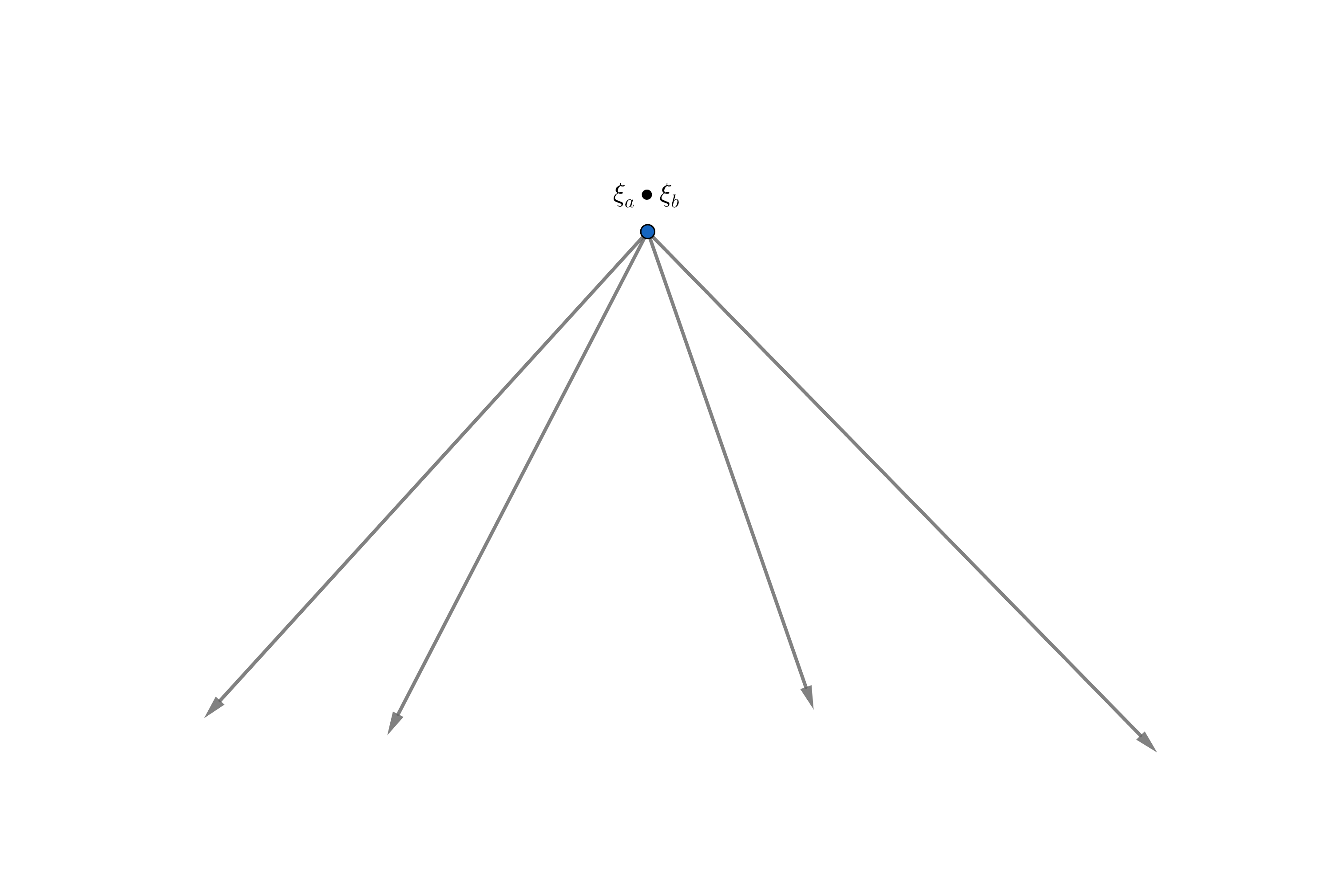} }}
    \caption{An exapmle of contracting edge between polyvector fields $\xi_a$ and $\xi_b$. Right hand side gives the graph after contraction}
    \label{fig:example}
\end{figure}
\newline
Summands in the second term of $(F)$ come from contraction of an admissible graph $\Gamma\in G_{n,\bar{n}}'$ by a subgraph $\Gamma_l$, subject to conditions:
\begin{enumerate}
    \item $\Gamma_l$ is admissible, has $l$ vertices of type $1$ and $\bar{l}$ vertices of type 2, and $2l+\bar{l}-2=\#E(\Gamma_l)$.
    \item The graph $\Gamma_k$ obtained by contracting $\Gamma_l$ to a new point of second type is admissible, has $k$ vertices of type $1$ and $\bar{k}$ vertices of type 2, and $2k+\bar{k}-2=\#E(\Gamma_k)$.
\end{enumerate}
\begin{figure}[htp]
    \centering
    \includegraphics[width=9cm]{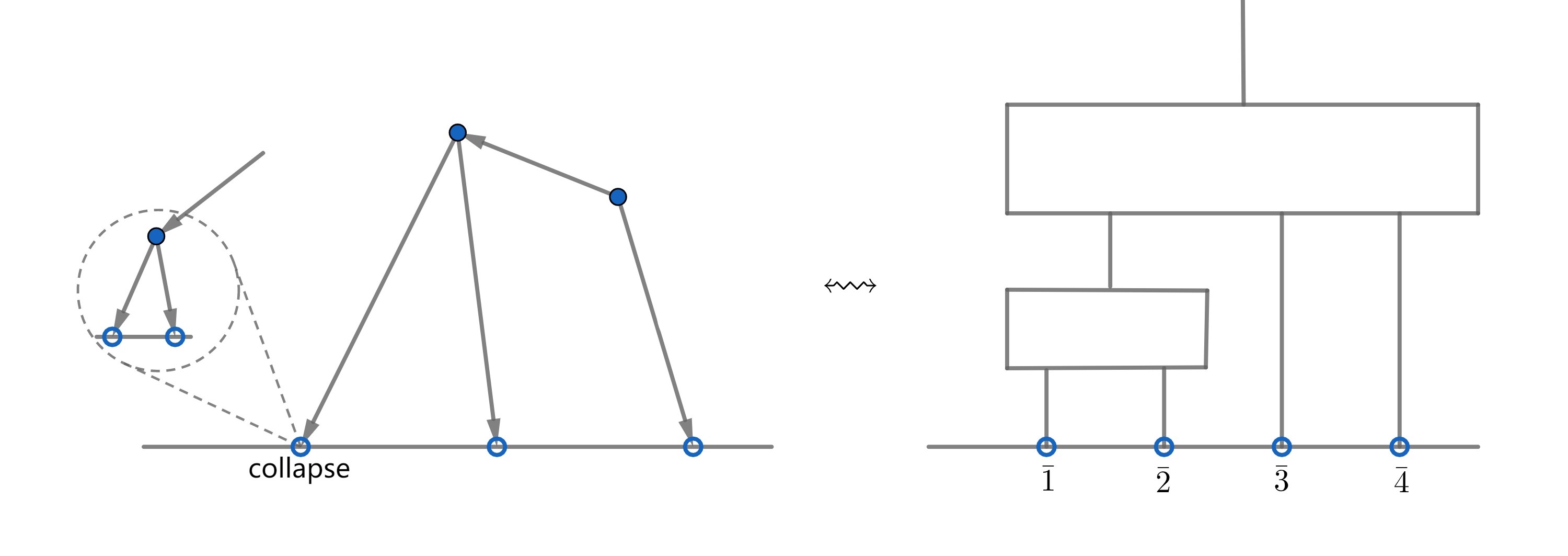}
    \caption{An example of a contraction corresponding to $\circ_{\bar{1}}$}
    \label{fig:Kontsevich's eye}
\end{figure}
Let $f:=(U_{\Gamma_k}(\xi_{\sigma(1)}\wedge\cdots\wedge\xi_{\sigma(k)})$, $g:=U_{\Gamma_l}(\xi_{\sigma(k+1)}\wedge\cdots\wedge\xi_{\sigma(n)}))$. The $\bar{i}^{th}$ circle product of $f$ and $g$ precisely gives the information of contracting $\Gamma_l$ to a vertex of second type in $\Gamma_k$; On the other hand a contraction of $\Gamma'$ satisfying $(1)(2)$ can be identified with decomposition of strata of type S2 (as discussed before), which can then be rearranged to a summand of the circle product.\\
Such $\Gamma'$ has $k+l$ vertices of first type, $\bar{k}+\bar{l}-1$ vertices of second type, and $\#E(\Gamma')=2(k+l)+(\bar{k}+\bar{l}-1)-3$.
\subsubsection{Verification of $(F)$.}
From previous discussion we can write the difference between LHS and RHS as \[\sum_{\Gamma'}c_{\Gamma'}U_{\Gamma'}(\xi_1\otimes\cdots\otimes\xi_n)(f_1\otimes\cdots\otimes f_{\bar{n}}).\]
for $\Gamma'\in G_{n,\bar{n}}$ with $\#E(\Gamma')=2n+\bar{n}-3$.
We are going to show that $C_{\Gamma'}=0$ for all such $\Gamma'$, the key is to use \textit{Stoke's theorem}.
\begin{theorem}
\textit{(Stoke's theorem on manifolds with corners)} Let $M$ be an oriented smooth manifold with corners of dimension $n$, and let $\omega$ be a compactly supported smooth $(n-1)$-form on $M$. Then
\[\int_Md\omega=\int_{\partial M}\omega.\]
\end{theorem}
For $\Gamma'\in G_{n,\bar{n}}$ such that $\#E(\Gamma')=2n+\bar{n}-3$, we define the $(2n+\bar{n}-3)$-form $\omega:=\bigwedge\limits_{e\in E(\Gamma')}d\phi_e$, where each $\phi_e$ is defined the same way as in $\S 6.2$. By Leibniz rule, $d\omega=0$. Meanwhile, since $\bar{C}^+_{n,\bar{n}}$ is compact and of dimension $2n+\bar{n}-2$, Stoke's theorem says
\[\int_{\partial \bar{C}^+_{n,\bar{n}}}\omega=\int_{\bar{C}^+_{n,\bar{n}}}d\omega,\]
which then vanishes as $d\omega=0$.\newline
We are going to decompose the boundary stratum $\partial \bar{C}^+_{n,\bar{n}}$ into strata of type S1 and S2. We then integrate $\omega$ separately on each stratum, and prove that total integration as in LHS above is equal to $c_{\Gamma'}$.

\textbf{The case S1.}
Type S1 strata can be described locally as $C_i\times C_{n-i+1,\bar{n}}$ for $i\geq 2$. Since the integral vanishes if the degree of differential form is not equal to the dimension of the underlying manifold, the contracted subgraph $\Gamma_1$ needs to have $dim(C_i)=2i-3$ edges.
We subdivide the situation into two cases:
\begin{enumerate}
    \item $i=2$: in this case there is only one edge $e$ in the first component. We are integrating $d\phi_e$ over $C_2$ which is homeomorphic to $S^1$. The integration gives a $\pm2\pi$ which then cancels out one of the $2\pi$ from the fraction of its associated $W_{\Gamma'}$ . The remaining form is given by a new graph $\Gamma$ obtained by contracting the edge $e$ in $\Gamma'$, which integrates (over $\Gamma=$ the second component) to give the weight $W_{\Gamma}$ of $U_{\Gamma}$ in the first term of $(F)$. Thus summing over all graphs and strata of this sort gives the first term of $(F)$.
    \item $i>2$: this results from the following lemma, which we only sketch the main idea at the end.
    The rigorous proof can be found in \cite{Kontsevich_2003}, see \cite{Khovanskii1997} for a proof in the algebro-geometric setting.
\end{enumerate}
\begin{lemma}
(\textit{The vanishing lemma}) Let $n\geq 3$ be an integer. The integral over $C_n$ of the wedge product of any $dim(C_n)=2n-3$ closed $1$-forms $d\phi_{i_{\alpha},j_{\alpha}},\alpha=1,\cdots,2n-3$, is equal to $0$.
\end{lemma}
\textbf{The case S2.}
This type of strata locally have the form
$C_{k,\bar{k}}\times C_{l,\bar{l}}$ where $k+l=n$, $\bar{l}=\bar{n}+1-\bar{k}$. The integral of the differential form associated to our graph decomposes into a product of two integrals. For dimension reason we require the subgraph $\Gamma_k\in G_{k,\bar{k}}$ has number of edges equal to $2k+\bar{k}-2$.\\
We call an edge starting from vertices of $\Gamma_k$ a \textit{bad edge} if it lands in a point not belong to $V_1(\Gamma_k)\sqcup V_2(\Gamma_k)$. The contracted graph of $\Gamma'$ by $\Gamma_k$ is admissible if and only if there is no bad edge occur in $\Gamma'$. \\
From analysis of $(F)$ we know that the contraction of $\Gamma'$ by $\Gamma_k$ contribute to the second term of $(F)$. If we sum over all possibilities of $k$ we obtain the second term. \\
\begin{figure}[htp]
    \centering
    \includegraphics[width=8cm]{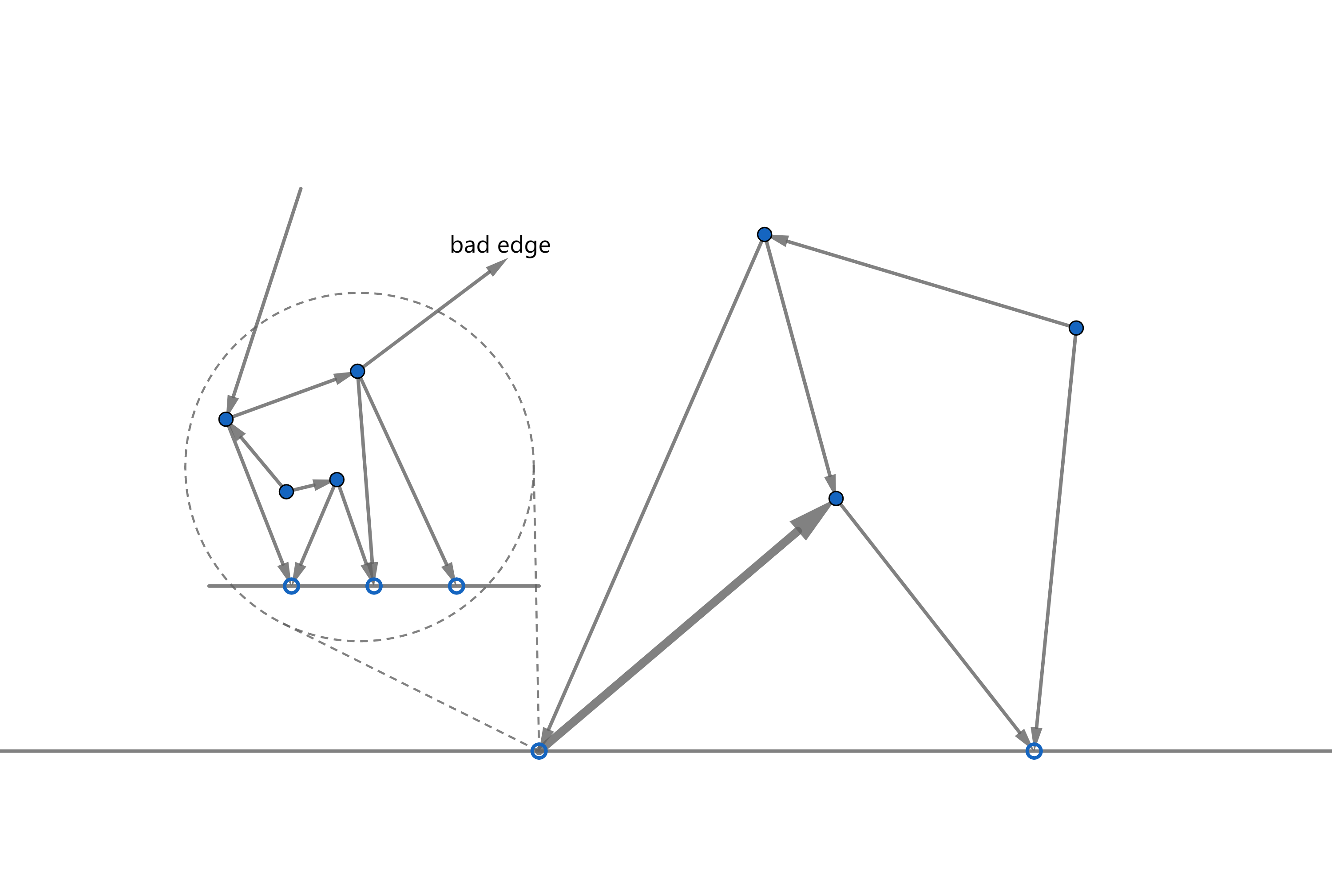}
    \caption{The resulting graph is not admissible if a bad edge occurs.}
    \label{Example of contraction with a bad edge}
\end{figure}\newline
If $\Gamma'$ has a bad edge, then $\Gamma_l$ is not admissible, for it has an edge $e$ start from a point of $2^{nd}$ type. But then $d\phi_e$ vanishes no matter how the endpoint of $e$ moves. Because, in our chosen metric system, geodesic lines starting from points on $x$-axis always have the angle measure 0.
\newline
We have exhausted all possible decomposition of $\partial \bar{C}^+_{n,\bar{n}}$. Once we proved the vanishing lemma, our proof of main theorem in $\mathbb{R}^d$ is complete.\newline\newline
\textit{Sketch of main idea for the vanishing lemma:} The main idea is, first, we restrict the integral to an even number of angle forms. This is done by identifying $C_n$ with with a subset of $(\mathbb{H}^2)^n$. Pick an edge $e_1$, use affine transformation to move one of its endpoints to the origin and the other to lie on the unit circle. Then the integral decomposes into a product of $d\phi_{e_1}$ integrated over $S^1$ (which gives $\pm 2\pi$) and the remaining $2N:=2n-4$ forms integrated over the resulting (complex) manifold $U$ given by the isomorphism $C_n\cong S^1\times U$.\newline
A problem occurs when integrating $\int_U\bigwedge\limits_{\alpha=1}^{2N}(dArgZ_{\alpha})$: $Arg(-)$ is not globally well-defined. To address the problem, recall that our angle map is harmonic, and from a complex analysis course we can find a harmonic conjugate $G(-)$ to $Arg(-)$ so that $Arg(-)+iG(-)$ is holomorphic. It turns out that the harmonic conjugate is $Log|-|$, which is globally well-defined. \\
We want to evaluate the integral of $\bigwedge\limits^{2N}_{\alpha=1}dLog|Z_{\alpha}|$ over $U$. Note that being a wedge product of $1$-forms, $\bigwedge\limits^{2N}_{\alpha=1}dLog|Z_{\alpha}|$ is exact and we can write it as $d\omega$ for some $2N-1$-form $\omega$.
At this stage, we denote $\bar{U}$ for the compactification of $U$, and integrate our form over $\bar{U}$:
\[\int_{\bar{U}}\bigwedge^{2N}_{\alpha=1}dLog|Z_{\alpha}|=\int_{\partial\bar{U}}\omega=0\]
Where the first equality comes from Stoke's theorem, and the second equality is by dimension reason. Finally, it was in Lemma 6.6 of \cite{Kontsevich_2003}
where Kontsevich show that the contribution of the integral over $\bar{U}-U$ is zero.
\hfill{\qed}
\begin{remark}
Using $(F)$, given any Poisson bivector field $\Pi$ on $\mathbb{R}^d$ we can write out the star product for $\Pi$ explicitly, as
\[P(\Pi):=\sum_{j\geq 0}\frac{\hbar^j}{j!}U_j(\Pi\wedge\cdots\wedge\Pi).\]
$P(\Pi)$ is associative: in $(F)$ take every $\xi$ to be $\Pi$. The first term of $(F)$ becomes zero while the second term reads $P(\Pi)(P(\Pi)\otimes id-id\otimes P(\Pi))=0$.
Expanding up to degree $1$ of $\hbar$, we have
\[f\otimes g\mapsto fg+\hbar\Pi(df,dg)+O(\hbar^2).\]
A more general fact behind: Any $L_{\infty}$ morphism $f:L\rightarrow L'$ an $L_{\infty}$ between DGLAs induces a morphism $MC(L\otimes m)\rightarrow MC(L'\otimes m)$, where $m$ is a nilpotent non-unital algebra, given by
\[s\mapsto \sum_{n=1}^{\infty}\frac{1}{n!}f_n(s\wedge\cdots\wedge s).\]
\end{remark}
\subsection{A view toward globalizing}
We have shown the morphism $U$ between $T_{poly}(\mathbb{R}^d)$ and $D_{poly}(\mathbb{R}^d)$ is $L_{\infty}$ quasi-isomorphic. The construction of these two DGLAs still work for subsets of $\mathbb{R}^d$. Taking formal completion $\mathbb{R}^d_{formal}$ of $\mathbb{R}^d$ at zero, $i.e.$ the formal spectrum of $\lim\limits_{\leftarrow}A/I^n$, where $I$ is the ideal of $A:=C^{\infty}(\mathbb{R}^d)$, we can similarly construct $T_{poly}(\mathbb{R}^d_{formal})$ and $D_{poly}(\mathbb{R}^d_{formal})$. Moreover, the morphism $U$ on $\mathbb{R}^d$ induces a morphism $\tilde{U}$ on $\mathbb{R}^d_{formal}$.
$\tilde{U}$ is $L_{\infty}$ quasi-isomorphic on $\mathbb{R}^d_{formal}$: $(F)$ is automatic for $\tilde{U}$, and HKR works for any manifold.\\
The point we are considering $\mathbb{R}^d_{formal}$ is that the completion of $\mathbb{R}^d$ at zero gives local data of $\mathbb{R}^d$. We want to glue these local data together to form the global data.\\
Let $W_d$ be the vector field over $\mathbb{R}^d_{formal}$. Regarding $W_d$ as a DGLA with trivial grading, we have natural morphism of DGLAs:
\[m_T:W_d\longrightarrow T_{poly}(\mathbb{R}^d_{formal})\]
viewing $W_d$ as the set of derivations, and
\[m_D:W_d\longrightarrow D_{poly}(\mathbb{R}^d_{formal})\]
viewing $W_d$ as a set of differential operators.\\
Let $G_d$ be Lie group of formal diffeomorphisms of $\mathbb{R}^d$ preserving base point $0$. If we can show that $\tilde{U}$ is \textit{compatible} with the actions of $W_d$ and $G_d$, then these local $L_{\infty}$ quasi-isomorphisms patch together to give a global $L_{\infty}$ quasi-isomorphism. Precisely, by compatible it means $\tilde{U}$ to satisfy:
\begin{enumerate}
    \item for any $\xi\in W_d$ we have equality
    \[\tilde{U}(m_T(\xi))=m_D(\tilde{U}(\xi));\]
    \item $\tilde{U}$ is $GL(d,\mathbb{R})$-equivariant (The action of $GL(d,\mathbb{R})$ on $\mathbb{R}^d_{formal}$ is induced by its action on $\mathbb{R}$);
    \item for any $k\geq 2$, $\xi_1,\cdots,\xi_k\in W_d$ we have the equality
    \[\tilde{U}_k(m_T(\xi_1)\otimes\cdots\otimes m_T(\xi_k))=0;\]
    \item for any $k\geq 2$, $\xi\in gl(d,\mathbb{R})\subset W_d$, and for any $\eta_2,\cdots,\eta_k\in T_{poly}(\mathbb{R}^d_{formal})$ we have
    \[\tilde{U}_k(m_T(\xi)\otimes\eta_2\otimes\cdots\otimes\eta_k)=0\]
\end{enumerate}
The detailed verification of $1--4$ listed above can be found in \cite{Kontsevich_2003} and \cite{Cattaneo_2001}.
\newpage
\appendices
\section{About bar complex}
We assume standard knowledge ((co)chain complex, (co)chain homotopy, free resolution, etc.) from a homological algebra course.\newline
$A$ is assumed to be an associative $k$-algebra.\\
As defined in $3.1$, a \textit{bar complex} is a sequence of $A$-bimodules $(A^{\otimes n})_{n\in\mathbb{N}}$ with $d_n:A^{\otimes n+2}\rightarrow A^{\otimes n+1} $ given by
\begin{equation*}
    a_0\otimes \cdots \otimes a_{n+1}\mapsto \sum^n_{i=0}(-1)^ia_0\otimes \cdots \otimes a_{i-1}\otimes a_ia_{i+1}\otimes a_{i+2}\otimes \cdots \otimes a_{n+1},
\end{equation*}
for $a_0,\cdots a_{n+1}\in A$.\\
Before showing that $(d_n)_{n\in \mathbb{N}}$ are differentials, we first introduce the morphisms:
\[r_n: A^{\otimes n+2}\rightarrow A^{\otimes n+3}: 1\otimes a_0\otimes\cdots\otimes a_{n+1}\]
It is clear that each $r_n$ is injective. Moreover, $Im(r_n)$ generates $A^{n+3}$ as a left A-module.
\[\begin{tikzcd}[row sep=1.5cm,column sep=1.5cm]
\cdots \arrow[r,"d_{n+1}"] & A^{\otimes n+2}\arrow[dl,"r_{n}" description]
\arrow[r,"d_n"] \arrow[d,"id"]
& A^{\otimes n+1} \arrow[dl,"r_{n-1}" description,swap]
\arrow[r,"d_{n-1}"] \arrow[d,"id"]
& A^{\otimes n} \arrow[dl,"r_{n-2}" description,swap]
\arrow[r,"d_{n-2}"] \arrow[d,"id"]
& \arrow[dl,"r_{n-3}" description,swap]\cdots\\
\cdots \arrow[r,"d_{n+1}"] &  A^{\otimes n+2}\arrow[r,"d_{n}"] & A^{\otimes n+1}
\arrow[r,"d_{n-1}"] & A^{\otimes n}
\arrow[r,"d_{n-2}"] & \cdots
\end{tikzcd}\]
A direct computation gives $d_{n+1}\circ r_n+r_{n-1}\circ d_n=id_{A^{\otimes n+2}}$:
\begin{align*}
    &\text{ }d_{n+1}\circ r_n(a_0\otimes\cdots\otimes a_n)+r_{n-1}\circ d_n(a_0\otimes\cdots\otimes a_n)\\
    &=a_0\otimes\cdots\otimes a_n+\sum_{i=1}^{n+1}(-1)^i1\otimes a_0\otimes\cdots\otimes a_n
    +\sum_{i=0}^{n}(-1)^i1\otimes a_0\otimes\cdots\otimes a_n\\
    &=a_0\otimes\cdots\otimes a_n.
\end{align*}
\begin{lemma}
$d_{n}\circ d_{n+1}=0$ for each $n\in\mathbb{N}$.
\begin{proof}
By induction. When $n=1$, for $a,b,c\in A$,
\[d_0\circ d_1(a,b,c)=d_0((a_0a_1)a_2-a_0(a_1a_2))=(a_0a_1)a_2-a_0(a_1a_2)=0,\]by associativity of $A$.\\
Suppose $d_{n-1}\circ d_{n}=0$. Then
\begin{align*}
d_n\circ d_{n+1}\circ r_n&=d_n\circ (id_{A^{\otimes n+2}}-r_{n-1}\circ d_n)\\
&=d_n-d_n\circ r_{n-1}\circ d_n\\
&=(id_{n+1}-d_n\circ r_{n-1})\circ d_n\\
&=r_{n-2}\circ d_{n-1}\circ d_n=0.
\end{align*}
Then $d_{n-1}\circ d_n=0$, since $Im(r_n)$ generates $A^{\otimes n+3}$.
\end{proof}
\end{lemma}
It follows that by taking $Hom_{A^e}(C^{bar}_{\bullet},M)$ or by tensoring an $A$-module $M$ respectively, the Hochschild (co)chain complexes are indeed (co)chain complexes.\\
Since $d\circ r+r\circ d=id=id-0$, $r_n$'s provides a chain homotopy between $id$ and the zero map. Since homotopic chain maps induces the same morphisms between (the chains of) homology groups, in particular, \[id_*=0_*:H_{\bullet}(C^{bar}(A))\rightarrow H_{\bullet}(C^{bar}(A)).\]
Hence $H_n(C^{bar}(A))$ is the $0$-module for each $n$, and the bar complex turns out to be exact. \\
\begin{lemma}
The component $C^{bar}_n= A^{\otimes n+2}$ of the chain complex is a free A-bimodule for any $n\geq 1$.
\begin{proof}
Take $\{a_i\text{ }|\text{ }i\in\Lambda\}$ as a vector space basis of $A^{\otimes n}$,
We have isomorphisms of $A$-bimodules:
\[A^{\otimes n+2}\cong A^e\otimes A^{\otimes n}\cong \bigoplus_{i\in\Lambda}A^e\cdot1\otimes a_i\otimes 1,\]
where the first isomorphism is given by
\[a_0\otimes\cdots\otimes a_{n+1}\mapsto (a_0\otimes a_{n+1})\otimes a_1\otimes\cdots\otimes a_n.\]
\end{proof}
\end{lemma}
\begin{theorem}
The bar complex $C^{bar}$ is a free resolution of $A$ as an $A$-bimodule, with augmentation $d_0:A\otimes A\rightarrow A$ the multiplication of $A$.
\begin{proof}
From above discussions we only need to show that the cokernel of $d_1:A\otimes A\otimes A\rightarrow A\otimes A$ is indeed the multiplication \[d_0:A\otimes A\rightarrow A: a_0\otimes a_1\mapsto a_0a_1.\]
We know that $d_0\circ d_1=0$, which induces a map
\[coker(d_1)\cong A\otimes A/Im(d_1)\rightarrow A: a_0\otimes a_1+Im(d_1)\mapsto a_0a_1,\]
whose inverse is given by
\[a\mapsto a\otimes 1:A\rightarrow coker(d_1).\]
\end{proof}
\end{theorem}
Thanks to the theorem, we can calculate left and right derived functors using bar complex. in particular, for an associative algebra $A$ and $A$-bimodule $M$,
\begin{align*}
    HH_n(A,M)\cong Tor_n^{A^e}(A,M).
\end{align*}
and
\begin{align*}
    HH^n(A,M)\cong Ext^n_{A^e}(A,M),
\end{align*}
This will be useful in our proof of HKR isomorphism.
\newpage
\section{A proof of HKR theorem}
We start by recalling some notions in Commutative algebra to set up the stage. Throughout the section, $A$ is assumed to be a noetherian commutative $k$-algebra. In this case, $A^{op}\cong A$ and $A^e\cong A\otimes_k A$. All the following terminology and results can be found in \cite{Eisenbud_1995} and \cite{weibel_1994}. \\
Recall that a \textit{Kähler differential} $\Omega_A$ is the $A$-module $I/I^2$, where $I$ is the kernel of the multiplication $\mu:A\otimes_k A\rightarrow A: a\otimes b\mapsto ab$. The Kähler differential co-represents the derivation, in the sense that $Hom_A(\Omega_A,M)\cong Der(A,M)$ for any $A$-module M.\newline
A \textit{regular sequence} in $A$ is a sequence $(a_1,\cdots,a_n)$ of elements in $A$ such that for each $i$, $a_i$ is not a zero divisor in $A/(a_1,\cdots,a_{i-1})$. A morphism of $A$-modules $f:M=A^{\oplus n}\rightarrow A$ is said to be \textit{associated to a regular sequence} $(a_1,\cdots,a_n)$ if the coefficients of $f$ are given by $(a_1,\cdots,a_n)$. In that case, the \textit{Koszul complex associated to $f$}, $i.e.$ the chain complex
\[
\begin{tikzcd}
0\overset{0}{\longrightarrow}\bigwedge\nolimits^{n}_{A}M\overset{d_f}{\longrightarrow} \bigwedge\nolimits^{n-1}_{A}M\overset{d_f}{\longrightarrow}\cdots\overset{d_f}{\longrightarrow}\bigwedge\nolimits^1_AM\overset{d_f}{\longrightarrow}A\longrightarrow 0
\end{tikzcd}
\]
whose the differential is given by
\[d_f:\bigwedge\nolimits^jM\rightarrow\bigwedge\nolimits^{j-1}M:m_1\wedge\cdots\wedge m_j\mapsto\sum_{i=1}^{j}f(m_i)m_1\wedge\cdots\wedge \hat{m_i}\wedge\cdots\wedge m_j \]
serves as a free resolution for $A/J$, where $J\lhd A$ is the ideal generated by the regular sequence $(a_1,\cdots,a_n)$.
\newline
Suppose $A$ is a flat $k$-algebra and (locally) of finite type over $k$. Then $A$ is \textit{smooth} over $k$ if, for all maximal ideals $m$ of $A$, the kernel of \[\mu_m:(A\otimes_k A)_{\mu^{-1}(m)}\longrightarrow A_m\] is generated by a regular sequence (where $\mu$ is the multiplication). In that case the Kähler differential $\Omega_A$ is a projective $A$-module of finite rank.
\newline
We now proof the central result of this section:
\begin{theorem}
(\textit{HKR theorem for the cochain complex}:) Suppose $A$ is a noetherian smooth commutative $k$-algebra, then there is a quasi-isomorphism \[\Phi:\bigwedge\nolimits^{\bullet}(Der(A))\rightarrow C^{\bullet}(A),\]
whose $0^{th}$ component is
$\Phi_0: Der(A)\rightarrow C^1(A): \xi\mapsto [a\mapsto \xi a]$
and whose $n^{th}$ component is
\begin{align*}
\Phi_{n}:\xi_0\wedge\cdots\wedge\xi_n\mapsto[(a_0\otimes\cdots\otimes a_n)\mapsto\frac{1}{(n+1)!}\sum_{\sigma\in S_{n+1}}sgn(\sigma)\prod_{i=0}^{n}\xi_{\sigma(i)}(a_i)].
\end{align*}
\begin{lemma}
Let $B$ be a commutative local ring, and $I\lhd B$ is an ideal generated by a regular sequence $g=(x_1,\cdots,x_n)$. Then the isomorphism
\begin{align*}
    \Phi_0: Hom_{B/I}(I/I^2,B/I)\overset{\cong}{\longrightarrow}Ext^1_B(B/I,B/I)
\end{align*}
can be extended into an isomorphism of graded commutative rings:
\begin{align*}
\Phi_{\bullet}:\bigwedge\nolimits^{\bullet}Hom_{B/I}(I/I^2,B/I)\overset{\cong}{\longrightarrow}Ext^{\bullet}_{B}(B/I,B/I).
\end{align*}
\begin{proof}
The Koszul complex provides a free resolution:
\[
\begin{tikzcd}
0\overset{0}{\longrightarrow}\bigwedge\nolimits^{n}_{B}B^{\oplus n}\overset{d_g}{\longrightarrow} \bigwedge\nolimits^{n-1}_{B}B^{\oplus n}\overset{d_g}{\longrightarrow}\cdots\overset{d_g}{\longrightarrow}\bigwedge\nolimits^1_BB^{\oplus n}\overset{d_g}{\longrightarrow}B\longrightarrow B/I\longrightarrow 0
\end{tikzcd}
\]
of $B/I$ as $B$-module, which we can use to compute $Ext$:
\begin{align*}
    Ext^{\bullet}_B(B/I,B/I)&\cong H^{\bullet}(Hom_B(\bigwedge\nolimits^{\bullet}_BB^{\oplus n},B/I),d:f\mapsto f\circ d_g)\\
    &\cong H^{\bullet}(Hom_B(\bigwedge\nolimits^{\bullet}_BB^{\oplus n},B/I),0)\\
    &\cong Hom_B(\bigwedge\nolimits^{\bullet}_BB^{\oplus n},B/I)\\
    &\cong \bigwedge\nolimits^{\bullet}_{B}Hom_B(B^{\oplus n},B/I)\\
    &\cong \bigwedge\nolimits^{\bullet}_{B/I}Hom_{B/I}((B/I)^{\oplus n},B/I)\\
    &\cong \bigwedge\nolimits^{\bullet}_{B/I}Hom_{B/I}(I/I^2,B/I).
\end{align*}
Where the $2^{nd}$ isomorphism is because all coefficients of $d_g$ land in $I$; the $4^{th}$ isomorphism is because $B^{\oplus n}$ is free of finite rank over $B$; the $5^{th}$ isomorphism is by changing point of view. The last isomorphism is because $I$ is generated by a regular sequence of length $n$, so $I/I^2$ is a free module over $B/I$ of rank $n$ (see \cite{matsumura_1987}).
Note that isomorphism in each step preserves grading.
\end{proof}
\end{lemma}
\begin{proof}\textit{(for the Theorem)}
For each $n$ we have constructed morphism \[\Phi_n:\bigwedge_A\nolimits^{n+1}Der(A)\rightarrow C^{n+1}(A)\]as the natural generalization of \[\Phi_0:Der(A)\rightarrow C^1(A).\] By local-global principle\footnote{Basically, to check that a homomorphism of $A$-modules $M\rightarrow N$ is isomorphic, one pass to check that it is isomorphic at every maximal ideal.} we can check that the induced morphism by $\Phi$ (still call it $\Phi$) on cohomology groups is isomorphic for each maximal ideal $m\lhd A$. That is, we want \[\Phi_n:\bigwedge\nolimits^{n+1}_ADer(A)\otimes_AA_m\rightarrow HH^{n+1}(A,A)\otimes_AA_m\]
is isomorphic for each maximal ideal $m\lhd A$.\\
We have for each $n\in \mathbb{N}$:
\begin{align*}
    LHS&=\bigwedge\nolimits_{A}^{n+1}Der(A)\otimes_AA_m\\
       &\cong\bigwedge\nolimits^{n+1}_{A}(Hom_A(\Omega_A,A))\otimes_A A_m\\
       &\cong Hom_A(\bigwedge\nolimits_{A}^{n+1}\Omega_A,A)\otimes_AA_m\\
       &\cong Hom_{A_m}(\bigwedge\nolimits^{n+1}_{A_m}\Omega_{A_m}, A_m)\\
       &\cong\bigwedge\nolimits^{n+1}_{A_m}Hom_{A_m}(\Omega_{A_m},A_m)\tag{L}
\end{align*}
Where the $1^{st}$ isomorphism comes from isomorphism $Hom_A(\Omega_A,A)\cong Der(A)$; the $2^{nd}$ is by the commutativity of taking wedge product and taking dual space; the $3^{rd}$ isomorphism is because $\Omega_A$ being a projective module of finite rank over $A$, it is locally free of finite rank, hence finitely presented, and so is $\bigwedge\nolimits ^{n+1}_{A}\Omega_{A}$. \\
If we can show that
\begin{align}
    RHS=HH^{n+1}(A,A)\otimes_AA_m\cong HH^{n+1}(A_m,A_m), \tag{R}
\end{align}
then by taking $B=A_m\otimes_k A_m$, $I=ker(\mu_m:A_m\otimes_k A_m\rightarrow A_m: \alpha\otimes \beta\mapsto \alpha\beta)$, we have
\begin{align*}
    HH^{n+1}(A_m)&\cong Ext^{n+1}_{A_m\otimes A_m}(B/I,B/I)\\
    &= Ext^{n+1}_B(B/I,B/I)\\
    &\cong \bigwedge\nolimits^{n+1}_{B/I}Hom_{B/I}(I/I^2,B/I)\\
    &\cong \bigwedge\nolimits^{n+1}_{B/I}Hom_{B/I}(\Omega_{A_m},A_m)\\
    &\cong (L).
\end{align*}
Where the $1^{st}$ isomorphism is just by plug-in, the $2^{nd}$ isomorphism is by above lemma, the $3^{rd}$ isomorphism comes from definition of $\Omega_A$.\\
To show $(R)$, notice it always holds that
\begin{align*}
HH^{n+1}(A,A)\otimes A_m&\cong Ext^{n+1}_{A\otimes_k A}(A,A)\otimes A_m\\
&\cong Ext^{n+1}_{A_m\otimes_k A_m}(A_m,A_m)\cong HH^{n+1}(A_m,A_m).
\end{align*}
So, indeed, $\Phi_0$ induces isomorphisms $(\bigwedge\nolimits_{A}^{n+1}Der(A)\rightarrow HH^{n+1}(A))_{n\in\mathbb{N}}$ in a natural way. Since $\Phi_n$ defined in Theorem 27 are exactly those natural extension of $\Phi_0$ to degree $n+1$, these $(\Phi_n)$ do give a quasi-isomorphism between $\bigwedge\nolimits^{\bullet}Der(A) $ and $C^{\bullet}(A)$.
\end{proof}
\end{theorem}
\newpage
\nocite{PBelmans}
\nocite{article}
\nocite{doubek2009deformation}
\nocite{fu2006algebraic}
\printbibliography
\end{document}